\documentclass[12pt]{amsart}
\usepackage{amsmath,amsfonts,amssymb,amscd,amsthm,amsbsy,epsf}
\newtheorem{thm}{Theorem}[section]
\newtheorem*{thmA}{Theorem A}
\newtheorem*{thmB}{Theorem B}

\newtheorem{lem}[thm]{Lemma}
\newtheorem{cor}[thm]{Corollary}
\newtheorem{prop}[thm]{Proposition}
\theoremstyle{definition} 		
\newtheorem{defn}[thm]{Definition}
\newtheorem{remark}[thm]{Remark}

\def\nat{{\mathbb N}}
\def\A{{\mathcal A}}
\def\B{{\mathcal B}}
\def\C{{\mathcal C}}

\def\F{{\mathcal F}}
\def\G{{\mathcal G}}

\def\R{{\mathcal R}}
\def\U{{\mathcal U}}

\def\W{{\mathcal W}}

\def\bs{{\mathbf s}}
\def\bt{{\mathbf t}}
\def\bu{{\mathbf u}}
\def\bw{{\mathbf w}}

\def\Tau{{\mathfrak T}}
\begin{document}
\title{Schreier Sets in Ramsey Theory}
\author{V. Farmaki and S. Negrepontis}

\begin{abstract}
We show that Ramsey theory, a domain presently conceived to guarantee the existence of large homogeneous sets for partitions on $k$-tuples of words (for every natural number $k$) over a finite alphabet, can be extended to one for partitions on Schreier-type sets of words (of every countable ordinal). Indeed, we establish an extension of the partition theorem of Carlson about words and of the (more general) partition theorem of Furstenberg-Katznelson about combinatorial subspaces of the set of words (generating from $k$-tuples of words for any fixed natural number $k$) into a partition theorem about combinatorial subspaces (generating from Schreier-type sets of words of order any fixed countable ordinal). Furthermore, as a result we obtain a strengthening of Carlson's infinitary Nash-Williams type (and Ellentuck type) partition theorem about infinite sequences of variable words into a theorem, in which an infinite sequence of variable words and a binary partition of all the finite sequences of words, one of whose components is, in addition, a tree, are assumed, concluding that all the Schreier-type finite reductions of an infinite reduction of the given sequence have a behavior determined by the Cantor-Bendixson ordinal index of the tree-component of the partition, falling in the tree-component above that index and in its complement below it.
\end{abstract}

\keywords{Ramsey theory, Schreier sets, words}
\subjclass{Primary 05D10; Secondary 28D05}

\maketitle
\baselineskip=18pt
\pagestyle{plain}              

\setcounter{section}{0}
\section*{Introduction}

Our aim is to extend Ramsey theory so that it applies not only to partitions of $k$-tuples of words but more generally to partitions of Schreier-type sets of words of a fixed countable ordinal number. For a finite non-empty alphabet $\Sigma$ we denote by $W^{k}(\Sigma)$ (respectively, $W^{k}(\Sigma ; \upsilon)$) the family of sequences of $k$ many words (respectively, variable words) over $\Sigma$, and by $W^{\omega}(\Sigma ; \upsilon)$ the family of infinite sequences of variable words over $\Sigma$. By a reduction (respectively, variable reduction) of $\vec{w}=(w_n)_{n\in\nat} \in W^\omega (\Sigma ; \upsilon)$ we mean any infinite sequence of words (respectively, variable words), denoted by $\vec{u}\prec\vec{w}$, obtained from $\vec{w}$ by replacing each occurence of the variable in each $w_n$ by one element of the set $\Sigma \cup \{\upsilon\}$, dividing the resulting sequence into infinitely many finite blocks of consecutive words, and concatenating the members of each block; the first element (respectively, the first $k$ elements) of a reduction of $\vec{w}$ is called a reduced word (respectively, a finite reductions with $k$ words) of $\vec{w}$. (These terms will be defined more formally below). For a natural number $r$, an $r$-coloring (or an $r$-partition) of a set $S$ is a map $\chi : S \rightarrow \{1,\ldots,r\}$, and $\chi(s)$ is the color of $s$ for $s\in S$. A set $T\subseteq S$ is monochromatic (under $\chi$) if $\chi$ is constant on $T$. 

The fundamental classical partition theorems of Ramsey theory, namely (a) Carlson's partition theorem (Lemma 5.9 in \cite{C}, Corollary 4.6 in \cite{BBH} in strenthened form), (b) the Furstenberg-Katznelson partition theorem (Theorems 2.7 and 3.1 in \cite{FK1}), and (c) Carlson's Nash-Williams type infinitary partition theorem (Theorem 2 in \cite{C}), can now be stated as follows:

\begin{thm}
[\textsf{Carlson's theorem}, \cite{C}, \cite{BBH}] 
\label{thm:ICarlson}
Let $\chi_1 : W^{1}(\Sigma)\rightarrow \{1,\ldots,r_1\}$ and $\chi_2 : W^{1}(\Sigma ; \upsilon)\rightarrow \{1,\ldots,r_2 \}$ be finite colorings of the sets $W^{1}(\Sigma)$ and $W^{1}(\Sigma ; \upsilon)$, respectively and $\vec{w} \in W^\omega (\Sigma ; \upsilon)$ be an infinite sequence of variable words over $\Sigma$; then there exists a variable reduction $\vec{u} \prec \vec{w} $ of $\vec{w}$ such that all the reduced words of $\vec{u}$ are  monochromatic under $\chi_1$ and all the reduced variable words of $\vec{u}$ are  monochromatic under $\chi_2$.
\end{thm}

\begin{thm}
[\textsf{Furstenberg-Katznelson's theorem}, \cite{FK1}]
\label{IFurst-Kat}
Let $k$ be any natural number, $\chi_1 : W^{k}(\Sigma)\rightarrow \{1,\ldots,r_1\}$ and $\chi_2 : W^{k}(\Sigma ; \upsilon) \rightarrow \{1,\ldots,r_2 \}$ be finite colorings of the sets $W^{k}(\Sigma)$ and $W^{k}(\Sigma ; \upsilon)$, respectively and $\vec{w} \in W^\omega (\Sigma ; \upsilon)$ be an infinite sequence of variable words over $\Sigma$; then there exists a variable reduction $\vec{u} \prec \vec{w} $ of $\vec{w}$ such that all the finite reductions with $k$ words of $\vec{u}$ are  monochromatic under $\chi_1$ and all the finite variable reductions with $k$ variable words of $\vec{u}$ are  monochromatic under $\chi_2$. 
\end{thm}
In addition Furstenberg and Katznelson in \cite{FK2} introduced the notion of a $k$- dimensional combinatorial subspace of $W(\Sigma)$ for $k$ any natural number and proved (in Theorem 3.1) a partition theorem about these combinatorial subspaces. 

\begin{thm} 
[\textsf{Carlson's infinitary partition theorem}, \cite{C}]
\label{cor:IblockNW}
Let  $\U \subseteq W^{\omega}(\Sigma ; \upsilon)$ be a pointwise closed family of infinite sequences of variable words over $\Sigma$ and $\vec{w} \in W^\omega (\Sigma ; \upsilon)$ be an infinite sequence of variable words over $\Sigma$;  then there exists a variable reduction $\vec{u}\prec\vec{w}$ of $\vec{w}$ over $\Sigma$ such that either
all the variable reductions of $\vec{u}$ are contained in $\U$ or all variable reductions of $\vec{u}$ are contained in the complement of $\U$.
\end{thm}

As stated, the aim of the present paper is to show that stronger versions of these partition theorems hold for the family of Schreier-type sets of words of every countable ordinal, and not just for the family of $k$-tuples of words, with $k$ restricted to a natural number. The hierarchy $(\A_\xi)_{\xi<\omega_1}$ of the families of Schreier sets of natural numbers, defined on the countable ordinals, provides a classification of the class of all finite subsets of the natural numbers measuring their complexity. The recursive definition of the Schreier sets $(\A_\xi)_{\xi<\omega_1}$ is as follows:  

\begin{defn}[\textsf{The Schreier system}, 
{[F1, Def. 7], [F2, Def. 1.5] [F3, Def. 1.3]}]
\label{Irecursivethin}
For every non-zero, countable, limit ordinal $\lambda$ choose and fix a strictly 
increasing sequence $(\lambda_n)_{n \in \nat}$ of successor ordinals smaller than $\lambda$ 
with $\sup_n \lambda_n = \lambda$. 
The system $(\A_\xi)_{\xi<\omega_1}$ is defined recursively as follows: 
\begin{itemize}
\item[(1)] $\A_0 = \{\emptyset\}$ and $\A_1 = \{\{n\} : n\in \nat\}$;
\item[(2)] $\A_{\zeta+1} = \{s\in [\nat]^{<\omega}_{>0} : s= \{n\} \cup s_1$, 
where $n\in \nat$, $\{n\} <s_1$ and $s_1\in \A_\zeta\}$; 
\item[(3i)] $\A_{\omega^{\beta+1}} = \{s\in [\nat]^{<\omega}_{>0} : 
s = \bigcup_{i=1}^n s_i$, where $n= \min s_1$, $s_1<\cdots < s_n$ and 
$s_1,\ldots, s_n\in \A_{\omega^\beta}\}$; 
\item[(3ii)] for a non-zero, countable limit ordinal $\lambda$, 
\newline
$\A_{\omega^\lambda} = \{s\in [\nat]^{<\omega}_{>0} : s\in \A_{\omega^{\lambda_n}}$
with $n= \min s\}$; and 
\item[(3iii)] for a limit ordinal $\xi$ such that $\omega^{\alpha}< \xi < \omega^{\alpha +1}$ for some $0< \alpha <\omega_1$, if \newline
$\xi = \omega^{\alpha} p
+ \sum_{i=1}^m \omega^{a_i} p_i$, where $m\in \nat$ with $m\ge0$, 
 $p,p_1,\ldots,p_m$ are natural numbers with $p,p_1,\ldots,p_m\ge1$ 
(so that either $p>1$, or $p=1$ and $m\ge 1$) and 
$a,a_1,\ldots,a_m$ are ordinals with $a>a_1>\cdots a_m >0$, 
\newline
$\A_\xi = \{s\in [\nat]^{<\omega}_{>0} :s= s_0 \cup (\bigcup_{i=1}^m s_i)$ 
with $s_m < \cdots < s_1 <s_0$,
$s_0= s_1^0\cup\cdots\cup s_p^0$ with $s_1^0<\cdots < s_p^0\in \A_{\omega^a}$,
and $s_i = s_1^i \cup\cdots\cup s_{p_i}^i$ with 
$s_1^i <\cdots<\ s_{p_i}^i\in \A_{\omega^{a_i}}$ $\forall\ 1\le i\le m\}$.
\end{itemize}
\end{defn} 
  
Note that in case 3(iii)) above the Cantor normal form of ordinals is employed (cf. \cite{KM}, \cite{L}).
Note also that for $k$ a natural number, i.e. a finite ordinal, the Schreier family $\A_k$ coincides with the family of all $k$-elements subsets of the natural numbers. 

In the definition of a Schreier system $(\A_\xi)_{\xi<\omega_1}$ we can fix for each non-zero, countable, limit ordinal $\lambda$ the  particular sequence $((\lambda)_n)_{n \in \nat}$ or equivalently the sequence of successor ordinals $((\lambda)_n^1)_{n \in \nat}$ defined below:

\begin{defn}\label{fixed}
Let $\lambda$ be a non-zero, countable, limit ordinal and ${n \in \nat}$.

{\rm (i)} We define inductively the ordinal $(\lambda)_n$ as follows: 
\begin{itemize}
\item[(1)] $(\omega)_n = n$. 
\item[(2)] $(\omega^{\alpha+1})_n = \omega^{\alpha}n$ for every $0 < \alpha< \omega_1$. 
\item[(3)] For a non-zero, countable limit ordinal $\alpha$ with $\alpha < \omega^{\alpha}$, 
\newline
$ (\omega^{\alpha})_n = \omega^{(\alpha)_n}$. 
\item[(4)] For a non-zero, countable limit ordinal $\alpha$ with $\alpha = \omega^{\alpha}$, 
let $\beta_0$ be the smallest ordinal such that $\alpha$ can be obtained as the limit of the sequence $(\beta_n)_{n \in \nat}$, where $\beta_{n+1} = \omega^{\beta_n}$. We set 
$(\alpha)_n = \beta_n$.
\item[(5)] For a limit ordinal $ \lambda = \omega^{\alpha_1}p_1 + \ldots + \omega^{\alpha_m}p_m $, where $m , p_1, \ldots , p_m \in\nat$ and  $ \lambda> \alpha_1 > \ldots > \alpha_m >0$ 
\newline
$ (\lambda)_n = \omega^{\alpha_1}p_1 + \ldots + \omega^{\alpha_{m-1}}p_{m-1} + 
\omega^{\alpha_m}(p_m-1) + (\omega^{\alpha_m})_n$. 
\end{itemize} 

{\rm (ii)}  Consider the strictly decreasing sequence $\lambda_{(0)}, \lambda_{(1)}, \ldots, \lambda_{(k_n^{\lambda})}$ defined by 
$\lambda_{(0)}= \lambda$, $\lambda_{(i+1)}=(\lambda_{(i)})_n$ for $0\leq i < k_n^{\lambda}-1$ and terminating when $\lambda_{(k_n^{\lambda})}$ is a successor ordinal. As this sequence of ordinals is strictly decreasing, it must terminate and let 
$ (\lambda)_n^1 = \lambda_{(k_n^{\lambda})}$ be the final term of this sequence. 
\end{defn}

Although the recursive Schreier system $(\A_\xi)_{\xi<\omega_1}$ is a purely 
combinatorial entity, it nevertheless arose gradually in 
connection with the theory of Banach  spaces, originally the family $\A_\xi$ was defined by Schreier 
(\cite{S}) (for $\xi=\omega$), next by Alspach-Odell \cite{AO} (for 
$\xi =\omega^\kappa$, $\kappa$ a natural number) and Alspach-Argyros
\cite{AA} (for $\xi=\omega^\alpha$, $\alpha$ a countable ordinal),
and finally by Farmaki \cite{F1}, \cite{F2}, \cite{F3} and Tomczak-Jaegermann \cite{TJ} 
(for $\xi$ any countable ordinal). (The reader is referred to the introduction of \cite{F3} for more 
details). 
\smallskip

Schreier sets were used firstly for the following transfinite extension of the classical Ramsey partition theorem (\cite {R}), a result about the existence of monochromatic sets for finite colorations of the family of all $k$-tuples, with $k$ a natural number:  

\begin{thm}
[\textsf{Ramsey partition theorem on Schreier sets}, (\cite{F2})]
\label{ordinal Ramsey}  
Let $\xi$ be a non-zero countable ordinal number. For any finite coloration $\chi$ 
of the family ${\A}_{\xi}$ and $M$ an infinite subset of $\nat$ there exists an infinite subset $L$ of $M$ such that ${\A}_{\xi} \cap [L]^{< \omega}$ is monochromatic.
\end{thm}

Using the family ${\A}_{\xi}$ we define (in Definition~\ref{recursivethin} ) the families $W^{\xi}(\Sigma)$,  $W^{\xi}(\Sigma ; \upsilon)$ of Schreier-type sets of words, variable words respectively over $\Sigma$, of a fixed countable ordinal number $\xi$. Carlson's theorem (Theorem~\ref{thm:ICarlson}) and the more general Furstenberg-Katznelson's theorem (Theorem~\ref{IFurst-Kat}) will be extended from $k$-tuples to Schreier-type sets of every countably ordinal; this is the content of the main Theorem in Section 2 (see Theorem~\ref{xi Carlson}). With the notation and definitions given in Section 1, it reads as follows:

\begin{thmA}
Let $\xi$ be a countable ordinal, $\chi_1 : W^{\xi}(\Sigma)\rightarrow \{1,\ldots,r_1\}$ and $\chi_2 : W^{\xi}(\Sigma ; \upsilon) \rightarrow \{1,\ldots,r_2 \}$ be finite colorings of the sets $W^{\xi}(\Sigma)$ and $W^{\xi}(\Sigma ; \upsilon)$, respectively and $\vec{w} \in W^\omega (\Sigma ; \upsilon)$ be an infinite sequence of variable words over $\Sigma$; then there exists a variable reduction $\vec{u} \prec \vec{w} $ of $\vec{w}$ such that 
all the finite reductions of $\vec{u}$ in the set $W^{\xi}(\Sigma)$ are monochromatic under $\chi_1$ and all the finite variable reductions of $\vec{u}$ in the set $W^{\xi}(\Sigma ; \upsilon)$ are monochromatic under $\chi_2$.
\end{thmA}

The proof of this result is closer to the method employed by us in proving Schreier type extensions of Hindman's and Milliken-Taylor's theorems in \cite{FN}, which in turn is inspired by the method invented by Baumgartner to prove Hindman's theorem in \cite{B}; in particular, we do not use topological dynamics (as employed in \cite{FK1}) or idempotent ultrafilters (as employed in \cite{C}, \cite{BBH}). Some consequences of the Main Theorem are  described in Section 2. Beside the Carlson and the Furstenberg-Katznelson theorems, Schreier-type extensions of the Hale-Jewett's theorem (\cite{HJ}) and consequently of the van der Waerden's theorem (\cite{vdW}) are obtained. 

Theorem A is next used, in conjuction with the tools developed in Section 3, one of which is a suitable Cantor- Bendixson index, to strengthen Carlson's infinitary theorem (Theorem~\ref{cor:IblockNW}) to various forms of Nash-Williams type partition theorems for words and variable words, involving Schreier families. A somewhat weaker version of our main results (Theorems~\ref{block-NashWilliams1}, ~\ref{block-NashWilliams2}) is contained in the following statement, which also strengthen Theorem~\ref{cor:IblockNW} (see Remark~\ref{rem:block-NashWilliams}).   

\begin{thmB}
Let $\G\subseteq W^{<\omega}(\Sigma)$ and $\F\subseteq W^{<\omega}(\Sigma ; \upsilon)$ be trees and $\vec{w} \in W^\omega (\Sigma ; \upsilon)$ be an infinite sequence of variable words over $\Sigma$; then either there exists a variable reduction $\vec{u} \prec \vec{w} $ of $\vec{w}$ such that all the finite reductions of $\vec{u}$ over $\Sigma$ are contained in $\G$ or there exists a countable ordinal $\xi_1 = \zeta_{\vec{w}}^{\G}$ 
such that for all $\xi >\xi_1$ there exists a variable reduction $\vec{u}\prec\vec{w}$ of $\vec{w}$ such that all the finite reductions of $\vec{u}$ in the set $W^{\xi}(\Sigma)$ are contained in the complement of $\G$. 

Furthermore, either there exists a variable reduction $\vec{u} \prec \vec{w} $ of $\vec{w}$ such that all the finite variable reductions of $\vec{u}$ over $\Sigma$ are contained in $\F$ or there exists a countable ordinal $\xi_2 = \zeta_{\vec{w}}^{\F}$ such that for all $\xi >\xi_2$ there exists a variable reduction $\vec{u}\prec\vec{w}$ of $\vec{w}$ such that all the finite variable reductions of $\vec{u}$ in the set $W^{\xi}(\Sigma ; \upsilon)$ are contained in the complement of $\F$.
\end{thmB}

Theorem B is strengthened, involving the Ellentuck topology $\Tau_E$ in Theorem~\ref{thm:Ellentuck}. A simple consequence of Theorem~\ref{thm:Ellentuck} is the characterization of 
completely Ramsey partitions of $W^{\omega}(\Sigma)$ and $W^{\omega}(\Sigma ; \upsilon)$ in terms of the Baire property 
in the topology $\Tau_E$, a result proved with different methods by Carlson in \cite{C}.

Let us remark at this point that the attractive alternative approach to Ramsey theory, via located words rather than 'classical' words, given by Bergelson-Blass-Hindman in \cite{BBH}, also admits of a Schreier-type extension, analogous to the one given in the present paper. The details will appear elsewhere.

The extended Ramsey theory developed in the present paper is a more powerful tool than the 'classical' Ramsey theory in that Schreier sets of all countable-ordinal orders capture a considerable part of Analysis, which is beyond the reach of the arithmetically oriented 'classical' Ramsey theory. This is attested by the fact that the Schreier families have found essential applications in Banach space theory on such questions as, for example, unconditionality, $l^1$ and $c_0$ embeddability, and distortion (see e.g. \cite{F1}, \cite{O}, \cite{AGR}, \cite{F4}). 

It is also noteworthy that the hereditary family $(\A_{\omega})_* = \{t\in [\nat]^{<\omega} : t\subseteq s$ 
for some $s\in \A_{\omega} \}\cup\{\emptyset\}$ generated by $\A_{\omega}$ figures prominently (under the name of the family of ``not large'' sets) in fundamental questions of mathematical logic, specifically in the (Ramsey type) 
Paris-Harrington statements, statements true and  provable in set-theory but unprovable in 
Peano arithmetic (cf. \cite{PH}, \cite{KS} and \cite{GRS}, pp. 169-180). The higher order hereditary Schreier families 
$(\A_{\xi})_*$ might  well be useful in forming and proving statements true but unprovable in certain systems 
endowed with induction stronger than that in Peano arithmetic.

The fact that Carlson's and Furstenberg-Katznelson's partition theorems have found important applications in various branches of Mathematics, including Ramsey ergodic theory, such as in the proof of the density version of the Hales-Jewett theorem by Furstenberg and Katznelson in \cite{FK2}, and the deep relations that exist between Ramsey ergodic theory and the exciting methods developed by Green and Tao in \cite{GT} on the existence of arbitrarily long arithmetic progressions of primes, make it reasonable to expect that the Schreier-type extension of Ramsey theory presented in this paper will find interesting applications.

\section{terminology and notation} 
We develop in this section the necessary terminology and notation.
We denote by $\nat=\{1,2,\ldots\}$ the set of natural numbers, $[\nat]^{<\omega}_{>0}$ the set of all non-empty, finite subsets of $\nat$, $[\nat]^{<\omega}=[\nat]^{<\omega}_{>0}\cup\{\emptyset\}$ and $[\nat]^{\omega}$ the set of all  infinite subsets of $\nat$.

Let $\Sigma$ be a finite, non-empty alphabet, and  $\upsilon\notin\Sigma$ an entity which we call a \textit{variable}. 
A \textit{word} over $\Sigma$ is a finite sequence of elements of $\Sigma$. The set of all the words over $\Sigma$ is denoted by $W(\Sigma)$ ; thus
\begin{center}
$W(\Sigma)= \{w=\alpha_1\ldots\alpha_k : k\in\nat,  \alpha_1,\ldots,\alpha_k\in \Sigma \}. $
\end{center}
$W(\Sigma)$ is turned into a semigroup by the operation of concatenation: the concatenation of two words $w_1=\alpha_1\ldots\alpha_k, w_2=\beta_1\ldots\beta_l$ over $\Sigma$ is defined to be the word 
\begin{center}
$w_1\ast w_2 = \alpha_1\ldots\alpha_k\beta_1\ldots\beta_l$.
\end{center}
For two words $w_1=\alpha_1\ldots\alpha_k, w_2=\beta_1\ldots\beta_l$ over $\Sigma$ we write 
\begin{center}
$w_1 \propto w_2$ iff $k<l$ and $\alpha_i= \beta_i$ for $i=1,\ldots,k$,  
\end{center}
and in case $w_1 \propto w_2$ we set $w_2-w_1=\beta_{k+1}\ldots\beta_l\in W(\Sigma)$.

A \textit{variable word} over $\Sigma$ is a word over $\Sigma\cup\{\upsilon\}$ in which $\upsilon$ actually appears. 
So, the set $W(\Sigma ; \upsilon)$ of variable words over $\Sigma$ is defined as
\begin{center}
$W(\Sigma ; \upsilon) = W(\Sigma\cup\{\upsilon\})\setminus W(\Sigma)$.
\end{center}
We note that the concatenation of two variable words is also a variable word. If $w$ is a variable word over $\Sigma$ and $\alpha \in\Sigma \cup\{\upsilon\}$ then we write $w(\alpha)$ for the result of replacing every occurence of the variable $\upsilon$ in $w$ by $\alpha$. Thus $w(\alpha)\in W(\Sigma)$ for $\alpha\in \Sigma$ and 
$w(\upsilon)=w$.
\newline
For two variable words $w_1=\alpha_1\ldots\alpha_k, w_2=\beta_1\ldots\beta_l$ over $\Sigma$ we write 
\begin{center}
$w_1 \propto w_2$ iff $k<l$, $\alpha_i= \beta_i$ for $i=1,\ldots,k$ and $w_2-w_1=\beta_{k+1}\ldots\beta_l\in W(\Sigma ; \upsilon)$. 
\end{center}

We denote by $W^{<\omega}(\Sigma )$ the family of all finite sequences of words over the alphabet $\Sigma$, by $W^{\omega}(\Sigma)$ the family of all infinite sequences of words over $\Sigma$ and by $W^{<\omega}(\Sigma ; \upsilon)$, $W^{\omega}(\Sigma ; \upsilon)$ the families of all finite, infinite 
sequences of variable words over $\Sigma$ respectively. Hence,
\begin{itemize}
\item[{}]$W^{<\omega}(\Sigma)= \{\bw=(w_1,\ldots,w_l): l\in\nat, w_1,\ldots,w_l\in W(\Sigma)\}\cup\{\emptyset\}$,
\item[{}]$W^{<\omega}(\Sigma ; \upsilon)= \{\bw=(w_1,\ldots,w_l): l\in\nat, w_1,\ldots,w_l\in W(\Sigma ; \upsilon)\}\cup\{\emptyset\}$,
\item[{}]$W^{\omega}(\Sigma)= \{\vec{w}=(w_n)_{n\in\nat}: w_n\in W(\Sigma)$ for every $n\in\nat \}$,
\item[{}]$W^{\omega}(\Sigma ; \upsilon)= \{\vec{w}=(w_n)_{n\in\nat}: w_n\in W(\Sigma ; \upsilon)$ for every $n\in\nat \}$.
\end{itemize}
The complexity of a finite sequence  $\bw=(w_1,\ldots,w_l)\in W^{<\omega}(\Sigma\cup\{\upsilon\})\setminus\{\emptyset\}$ of words, with $w_i=\alpha_{k_i}\alpha_{k_i+1}\ldots\alpha_{k_{i+1}-1}$ for $i=1,\ldots,l$, 
is described by the complexity of the corresponding finite sequence of natural numbers $1=k_1<\cdots<k_l<k_{l+1}\in\nat$, a complexity that will be described by the Schreier hierarchy; we thus define the correspondence 
\begin{center}
$d : W^{<\omega}(\Sigma\cup \upsilon)\setminus\{\emptyset\} \rightarrow [\nat]^{<\omega}$ such that $\bw=(w_1,\ldots,w_l)\rightarrow d(\bw)$
\end{center}
\begin{center}
with 
$d(\bw)= \emptyset$ if $l=1$, and 
$d(\bw) = \{k_2 < k_3 < \cdots < k_l \}$ if $l>1$.
\end{center}

Analogously, for every infinite sequence $\vec{w}=(w_n)_{n\in\nat}\in W^{\omega}(\Sigma\cup\{\upsilon\})$ of words, with $w_n=\alpha_{k_n}\alpha_{k_n+1}\ldots\alpha_{k_{n+1}-1}$ for all $n\in \nat$, the corresponding complexity is described by the complexity of the infinite sequence $1=k_1<k_2<k_3<\cdots\in \nat$ of natural numbers; we thus define the correspondence 
\begin{center}
$d : W^{\omega}(\Sigma\cup \upsilon)\rightarrow [\nat]^{\omega}$ 
with $d((w_n)_{n\in\nat}) = (k_2 <k_3 < \cdots < k_n < \cdots)$.
\end{center} 

A finite sequence $\bw=(w_1,\ldots,w_l)\in W^{<\omega}(\Sigma\cup\{\upsilon\})$ is an initial segment of the finite sequence $\bu=(u_1,\ldots,u_k)\in  W^{<\omega}(\Sigma\cup\{\upsilon\})$ iff $l<k$ and $u_i= w_i$ for $i=1,\ldots,l$ and $\bw$ is an initial segment of the infinite sequence $\vec{u}=(u_n)_{n\in\nat}\in W^{\omega}(\Sigma\cup\{\upsilon\})$ if $u_i= w_i$ for all $i=1,\ldots,l$. In these cases, extending to sequences the previous notation for words, we write $\bw\propto \bu$ or $\bw\propto \vec{u}$, and we set
\begin{center}
$\bu\setminus \bw = (u_{l+1},\ldots,u_k)$ and $\vec{u}\setminus \bw = (u_n)_{n>l}$ respectively.
\end{center}

\begin{defn} 
\label{reduction}
(1) (\textsf{Reduction of an infinite sequence of words by a word}) 
For an infinite sequence $\vec{w}=(w_n)_{n\in\nat} \in W^\omega (\Sigma ; \upsilon)$ of variable words and for a (variable or non-variable) word $t=\alpha_1\ldots \alpha_k\in\W(\Sigma\cup \{\upsilon\})$ (over the alphabet $\Sigma$), we set 
\begin{center}
$\vec{w}[t]=w_1(\alpha_1)\ast \ldots \ast w_k(\alpha_k)\in\W(\Sigma\cup \{\upsilon\})$.
\end{center}
The family $RW(\vec{w})$ of all \textit{the reduced words} and the family $VRW(\vec{w})$ of all the \textit{variable reduced words}  of $\vec{w}$ over $\Sigma$ are defined as follows:
\begin{center}
$RW(\vec{w})=\{\vec{w}[t] : t\in W(\Sigma)\}$ and $VRW(\vec{w})=\{\vec{w}[t] : t\in W(\Sigma ; \upsilon)\}$.
\end{center}
For $u_1=\vec{w}[t_1], u_2=\vec{w}[t_2] \in RW(\vec{w})$ (resp. $u_1,u_2 \in VRW(\vec{w}) )$ we write 
\begin{center}
$u_1 \propto u_2 $ iff $t_1 \propto t_2$.
\end{center}

(2) (\textsf{Reduction of an infinite sequence of words by a finite sequence of words}) 
For an infinite sequence $\vec{w}=(w_n)_{n\in\nat} \in W^\omega (\Sigma ; \upsilon)$ of variable words and for a finite sequence of (variable or non-variable) words $\bt=(t_1,\ldots,t_l)\in W^{<\omega}(\Sigma\cup\{\upsilon\})\setminus\{\emptyset\}$, with $t_i =\alpha_{k_i}\alpha_{k_i+1}\ldots\alpha_{k_{i+1}-1}$ for all $i=1,\ldots,l$, we set 
\begin{center}
$\vec{w}[\bt]=(u_1,\ldots,u_l)\in W^{<\omega}(\Sigma\cup\{\upsilon\})$ where $u_i = w_{k_i}(\alpha_{k_i})\ast w_{k_i+1}(\alpha_{k_i+1})\ast \ldots \ast w_{k_{i+1}-1}(\alpha_{k_{i+1}-1})$ for all 
$i=1,\ldots,l$.
\end{center}
Also we set $\vec{w}[\emptyset]=\emptyset$. The finite sequences of words $\vec{w}[\bt]$ for $\bt\in  W^{<\omega}(\Sigma)$ are called \textit{finite reductions } of $\vec{w}$ over $\Sigma$ and the finite sequences of variable words $\vec{w}[\bt]$ for $\bt\in W^{<\omega}(\Sigma ; \upsilon)$ are called \textit{finite variable reductions } of $\vec{w}$ over $\Sigma$. The set of all the finite reductions and the set of all the finite variable reductions of $\vec{w}$ over $\Sigma$ are denoted as follows:
\begin{center}
$RW^{<\omega}(\vec{w})=\{\vec{w}[\bt] : \bt\in W^{<\omega}(\Sigma)\}$ and $VRW^{<\omega}(\vec{w})=\{\vec{w}[\bt] : \bt\in W^{<\omega}(\Sigma ; \upsilon)\}$.
\end{center}
We set
\begin{center}
$d_{\vec{w}} : RW^{<\omega}(\vec{w})\cup VRW^{<\omega}(\vec{w})\setminus\{\emptyset\}\rightarrow [\nat]^{<\omega}$ 
with $d_{\vec{w}}(\vec{w}[\bt]) = d(\bt)$.
\end{center}
Observe that $RW^{<\omega}(\vec{e})= W^{<\omega} (\Sigma )$, $VRW^{<\omega}(\vec{e})= W^{<\omega} (\Sigma ; \upsilon)$ and $d_{\vec{e}}=d$ if $\vec{e}=(e_n)_{n\in\nat}$ with $e_n=\upsilon$ for every $n\in \nat$. Note also that it is not always true that $d_{\vec{w}}(\bu)= d(\bu)$ for every $\bu\in RW^{<\omega}(\vec{w})$.

(3) (\textsf{Reduction of an infinite sequence of words by an infinite sequence of words}) 
For an infinite sequence $\vec{w}=(w_n)_{n\in\nat} \in W^\omega (\Sigma ; \upsilon)$ of variable words and for
an infinite sequence 
$\vec{t}=(t_n)_{n\in\nat} \in W^\omega (\Sigma\cup\{\upsilon\})$ of (variable or non-variable) words, with 
$t_n=\alpha_{k_n}\alpha_{k_n+1}\ldots\alpha_{k_{n+1}-1}$ for all $n\in \nat$, we set
\begin{center}
$\vec{w}[\vec{t}]=(u_n)_{n\in\nat}\in \W^{\omega}(\Sigma\cup\{\upsilon\})$ where $u_n=w_{k_n}(\alpha_{k_n})\ast w_{k_n+1}(\alpha_{k_n+1})\ast \ldots \ast w_{k_{n+1}-1}(\alpha_{k_{n+1}-1})$ for all 
$n\in \nat$.
\end{center}
An infinite sequences of words $\vec{w}[\vec{t}]$ for $\vec{t}\in W^{\omega}(\Sigma)$ is called a \textit{reduction } of $\vec{w}$ over $\Sigma$ and for $\vec{t}\in  W^{\omega}(\Sigma ; \upsilon)$ a \textit{variable reduction } of $\vec{w}$ over $\Sigma$, respectively. The sets of all the reductions and all the variable reductions of $\vec{w}$ over $\Sigma$ respectively are denoted as follows:
\begin{center}
$RW^{\omega}(\vec{w})=\{\vec{w}[\vec{t}] : \vec{t}\in W^{\omega}(\Sigma)\}$ and $VRW^{\omega}(\vec{w})=\{\vec{w}[\vec{t}] : \vec{t}\in W^{\omega}(\Sigma ; \upsilon)\}$.
\end{center}
For $\vec{u}, \vec{w} \in W^\omega (\Sigma ; \upsilon)$, we write
\begin{center}
$\vec{u} \prec \vec{w}$ if and only if $\vec{u} \in VRW^{\omega}(\vec{w})$.
\end{center}
Notice that $\vec{u} \prec \vec{w}$ if and only if $VRW(\vec{u})\subseteq VRW(\vec{w})$. Hence, $\vec{w} \prec \vec{e}$ for every $\vec{w}\in W^{<\omega} (\Sigma ; \upsilon)$ in case 
$\vec{e}=(e_n)_{n\in\nat}$ with $e_n=\upsilon$ for every $n\in \nat$. We define 
\begin{center}
$d_{\vec{w}} : RW^{\omega}(\vec{w})\cup VRW^{\omega}(\vec{w})\rightarrow [\nat]^{\omega}$ 
with $d_{\vec{w}}(\vec{w}[\vec{t}]) = d(\vec{t})$.
\end{center}

(1*) (\textsf{Reduction of a finite sequence of words by a word}) For a finite sequence $\bw=(w_1,\ldots,w_n)\in \W^{<\omega} (\Sigma ; \upsilon)$ of variable words over the alphabet $\Sigma$, we define the sets
\begin{center}
$RW((w_1,\ldots,w_n))=\{w_1(\alpha_1)\ast \ldots \ast w_n(\alpha_n): \alpha_1\ldots\alpha_n \in W(\Sigma) \}$ and
\end{center}
\begin{center}
$VRW((w_1,\ldots,w_n))=\{w_1(\alpha_1)\ast \ldots \ast w_n(\alpha_n): \alpha_1\ldots\alpha_n \in W(\Sigma ; \upsilon) \}$,
\end{center}
of all the \textit{reduced words} and \textit{variable reduced words}, respectively, of $(w_1,\ldots,w_n)$ over $\Sigma$.

Notice that $RW((w_1,\ldots,w_n))$, $VRW((w_1,\ldots,w_n))$ are finite sets and that for a sequence $\vec{w}=(w_n)_{n\in\nat} \in W^\omega (\Sigma ; \upsilon)$ we have that $RW(\vec{w})=\bigcup\{RW(w_1,\ldots,w_n) : n\in \nat \}$ and $VRW(\vec{w})=\bigcup\{VRW(w_1,\ldots,w_n) : n\in \nat \}$.

(2*) (\textsf{Reduction of a finite sequence of words by a finite sequence of words}) For a finite sequence $\bw=(w_1,\ldots,w_n)\in W^{<\omega} (\Sigma ; \upsilon)$ of variable words over the alphabet $\Sigma$, we define analogously the families $RW^{<\omega}(\bw)$ and $VRW^{<\omega}(\bw)$ of all \textit{finite reductions} and \textit{variable finite reductions}, respectively, of $\bw$ over $\Sigma$. So, $(u_1,\ldots,u_l)\in RW^{<\omega}(\bw)$  if there exists $\bt=(t_1,\ldots,t_l)\in W^{<\omega}(\Sigma)$, where  $t_i=\alpha_{k_i}\alpha_{k_i+1}\ldots\alpha_{k_{i+1}-1}$ for all $i=1,\ldots,l$ and $k_{l+1} = n+1$, such that 
\begin{center}
$u_i = w_{k_i}(\alpha_{k_i})\ast w_{k_i+1}(\alpha_{k_i+1})\ast \ldots \ast w_{k_{i+1}-1}(\alpha_{k_{i+1}-1})$ for 
$i=1,\ldots,l$.
\end{center}
We set $d_{\bw}(\bu)= \{k_2,\ldots, k_{l}\}$.
\end{defn}
 
In the sequel we will also employ the following notation. For the families $\G\subseteq W^{<\omega} (\Sigma)$, $\F\subseteq W^{<\omega} (\Sigma ; \upsilon)$ and the words $s\in W(\Sigma)$, $t\in W(\Sigma ; \upsilon)$ we set
\begin{center}
$\G(s) = \{\bw\in W^{<\omega} (\Sigma) : $ either $\bw=(w_1,\ldots,w_l)\neq \emptyset$, $s \propto w_1$ and 
\end{center}
\begin{center}
$(s,w_1-s,w_2, \ldots,w_l) \in \G$ or $\bw=\emptyset$ and $(s) \in\G\}$, and
\end{center}

\begin{center}
$\F(t) = \{\bw\in W^{<\omega} (\Sigma ; \upsilon) : $ either $\bw=(w_1,\ldots,w_l)$, $t \propto w_1$ and  
\end{center}
\begin{center}
$(t,w_1-t,w_2, \ldots,w_l) \in \F$, or $\bw=\emptyset$ and $(t) \in\F\}$.
\end{center}
Also,
\begin{center}
$\G -s = \{\bw \in \G :$ either $\bw = (w_1,\ldots,w_l)$ and $s \propto w_1$, or 
$\bw =\emptyset\}$ and
\end{center}
\begin{center}
$\F -t = \{\bw \in \F :$ either $\bw = (w_1,\ldots,w_l)$ and $t \propto w_1$, or 
$\bw =\emptyset\}$.
\end{center}

For the sequence $\vec{w}=(w_n)_{n\in\nat} \in W^\omega (\Sigma ; \upsilon)$ and the words $t\in VRW(\vec{w})$,  $s\in RW(\vec{w})$) with $t\in VRW((w_1,\ldots,w_k))$ and $s\in RW((w_1,\ldots,w_k))$ for some $k \in \nat$, we set  
\begin{center}
$\vec{w} -t = (w_0, w_{k+2},w_{k+3},\ldots)\in VRW^{\omega}(\vec{w})$, where $w_0= t \ast w_{k+1}$, and 
$\vec{w} -s = (w_0, w_{k+2},w_{k+3},\ldots)\in VRW^{\omega}(\vec{w})$, where $w_0= s \ast w_{k+1}$.
\end{center}
Also, we set $\vec{w}\setminus t = \vec{w}\setminus s = (w_{k+1},w_{k+2},\ldots)\in \W^\omega (\Sigma ;\upsilon)$.

\section{The main partition theorem on Schreier families}

The main theorem of this section is Theorem~\ref{xi Carlson}, given in equivalent form in Theorem~\ref{xi-Furst-Kat}; this is a partition theorem for the Schreier finite sequences of words and the Schreier finite sequences of variable words over a finite non-empty alphabet $\Sigma$ of every countable order, and constitutes an extension to every countable order $\xi$ (a) of Carlson's theorem, Theorem~\ref{thm:ICarlson}, corresponding to ordinal level $\xi=0$ and (b) of Theorem~\ref{IFurst-Kat}, proved by Furstenberg and Katznelson, corresponding to finite ordinals $\xi < \omega$. 

In order to state Theorem~\ref{xi Carlson} we need the following definitions:

\begin{defn} 
[\textsf{The Schreier systems} $(W^{\xi}(\Sigma))_{0\leq \xi<\omega_1}$ and 
$(W^{\xi}(\Sigma ; \upsilon))_{0\leq \xi<\omega_1}$] 
\label{recursivethin}
Let $(\A_\xi)_{\xi<\omega_1}$ be a Schreier system of families of finite 
subsets of $\nat$ and $\Sigma$ be an alphabet. We will define the families $W^{\xi}(\Sigma)$ and $W^{\xi}(\Sigma ; \upsilon)$ of the Schreier finite sequences of words and of variable words over $\Sigma$ respectively, for every countable ordinal $\xi$ recursively as follows:
\begin{center}
$W^0(\Sigma)= \{\bw=(w_1) : w_1 \in W(\Sigma)\}$ and $W^0(\Sigma ; \upsilon)= \{\bw=(w_1) : w_1 \in W(\Sigma ; \upsilon)\}$,
\end{center}
and for every countable ordinal $\xi\ge1$ 
\begin{center}
$W^{\xi}(\Sigma) = \{\bw=(w_1,\ldots,w_l)\in W^{<\omega} (\Sigma) : d((w_1,\ldots,w_l))\in \A_\xi \}$; and
\end{center}
\begin{center}
$W^{\xi}(\Sigma ; \upsilon) = \{\bw=(w_1,\ldots,w_l)\in W^{<\omega} (\Sigma ; \upsilon) : d((w_1,\ldots,w_l))\in \A_\xi \}$.
\end{center}

For an infinite sequence $\vec{w}=(w_n)_{n\in\nat} \in \W^\omega (\Sigma ; \upsilon)$ of variable words over $\Sigma$,  we define the families of Schreier finite reductions and of variable reductions of $\vec{w}$ over $\Sigma$ as follows: 
\begin{center}
$RW^0(\vec{w})= \{\bu=(u_1) : u_1\in RW(\vec{w})\}$ and $VRW^0(\vec{w})= \{\bu=(u_1) : u_1\in VRW(\vec{w})\}$.
\end{center}
and for every countable ordinal $\xi\ge1$ 
\begin{center}
$RW^{\xi}(\vec{w}) = \{\bu=(u_1,\ldots,u_l)\in RW^{<\omega}(\vec{w}) : d_{\vec{w}}((u_1,\ldots,u_l))\in \A_\xi \}$; and
\end{center}
\begin{center}
$VRW^{\xi}(\vec{w}) = \{\bu=(u_1,\ldots,u_l)\in VRW^{<\omega}(\vec{w}) : d_{\vec{w}}((u_1,\ldots,u_l))\in \A_\xi \}$.
\end{center}
\end{defn}

Hence, a finite sequence $\bw\in W^{<\omega}(\Sigma)$ of words over $\Sigma$ belongs to the family $W^{\xi}(\Sigma)$ for some $1\leq\xi<\omega_1$ iff $\bw=(w_1,\ldots,w_l)$ with $l>1$ and there exist  $1=k_1<\cdots<k_l<k_{l+1}\in\nat$ with $\{k_i: 2\leq i\leq l\}\in \A_\xi$ and $ \alpha_1,\ldots,\alpha_{k_{l+1}-1}\in \Sigma$ such that 
$w_i =\alpha_{k_i}\alpha_{k_i+1}\ldots\alpha_{k_{i+1}-1}$ for all $i=1,\ldots,l$. 
  
Observe that $\bw=(w)\in W^0(\Sigma)$ for every $w\in W(\Sigma)$, while $\bw=(w)\notin W^{\xi}(\Sigma)$ for every  $\xi>0$, also that $W^{\xi}(\Sigma) = RW^{\xi}(\vec{e})$ and 
$W^{\xi}(\Sigma ; \upsilon) = VRW^{\xi}(\vec{e})$ for every countable ordinal $\xi$, in case  $\vec{e}=(e_n)_{n\in\nat}$ with $e_n=\upsilon$ for every $n\in \nat$ and that it is not true that 
$W^{\xi}(\Sigma) = RW^{\xi}(\vec{w})$ for every $\vec{w}\in W^\omega (\Sigma ; \upsilon)$.

The following proposition justifies the recursiveness of the Schreier systems $(W^{\xi}(\Sigma))_{0\leq \xi<\omega_1}$ and $(W^{\xi}(\Sigma ; \upsilon))_{0\leq \xi<\omega_1}$.

\begin{prop}\label{justification 2}
For every countable ordinal $\xi> 0$ there exists a concrete sequence $(\xi_n)_{n>1}$
of countable ordinals with $\xi_n<\xi$ for every $n\in \nat$, $1<n$ such that
\begin{center}
$W^{\xi}(\Sigma)(s) = W^{\xi_n}(\Sigma)\cap (W^{<\omega}(\Sigma) - s)$ and 
$W^{\xi}(\Sigma ; \upsilon)(t) = W^{\xi_n}(\Sigma ; \upsilon)\cap (W^{<\omega}(\Sigma ; \upsilon) - t)$
\end{center}
for every $n\in \nat$, $1<n$ and 
$s=\alpha_1\ldots\alpha_{n-1} \in W(\Sigma)$, $t=\beta_1\ldots\beta_{n-1} \in W(\Sigma ; \upsilon)$.

Moreover, $\xi_n =\zeta$ for every $n\in\nat$ in case $\xi = \zeta+1$ and 
$(\xi_n)$ is a strictly increasing sequence with $\sup_n \xi_n=\xi$ in case $\xi$ 
is a limit ordinal.
\end{prop} 
\begin{proof}
According to Proposition~1.7 in \cite{F3}, for every countable ordinal $\xi >0$ there exists a concrete sequence 
$(\xi_n)$ of countable ordinals with  $\xi_n <\xi$ such that 
\begin{equation*}
\A_\xi (n) = \A_{\xi_n} \cap [\{n+1,n+2,\ldots\}]^{<\omega}\ \text{ for every }
\ n\in \nat\ ,
\end{equation*}
where, $\A_\xi (n) = \{s \in [\nat]^{<\omega} : \{n\}<s ,  \{n\}\cup s \in \A_\xi \}$ for every $n\in \nat$.
\newline
Moreover, $\xi_n = \zeta$ for every $n\in\nat$ if $\xi = \zeta+1$ and 
$(\xi_n)$ is a strictly increasing sequence with $\sup_n \xi_n=\xi$ if $\xi$ 
is a limit ordinal. 

Let $n>1$ and $s=\alpha_1\ldots\alpha_{n-1} \in W(\Sigma)$, $t=\beta_1\ldots\beta_{n-1} \in W(\Sigma ; \upsilon)$.  We will prove that $W^{\xi}(\Sigma)(s) = W^{\xi_n}(\Sigma)\cap (W^{<\omega}(\Sigma) - s)$ for every countable ordinal $\xi> 0$. Similarly can be proved that $W^{\xi}(\Sigma ; \upsilon)(t) = W^{\xi_n}(\Sigma ; \upsilon)\cap (W^{<\omega}(\Sigma ; \upsilon) - t)$ for every $0<\xi<\omega_1$.  

For $\xi=1$, of course $W^{1}(\Sigma)(s) = W^{0}(\Sigma)\cap (W^{<\omega}(\Sigma) - s)$. Let $1<\xi<\omega_1$. Then 
$\emptyset \notin \W^{\xi}(\Sigma)(s)$, since if $\emptyset \in \W^{\xi}(\Sigma)(s)$, then $(s)\in W^{\xi}(\Sigma)$ for $\xi>0$ and of course $\emptyset\notin W^{\xi_n}(\Sigma)$ for every $0\leq \xi_n<\omega_1$. Let $\bw=(w_1,\ldots,w_l)\in W^{<\omega}(\Sigma)\setminus\{\emptyset\}$. Then there exist $1=k_1<\cdots<k_l<k_{l+1}\in\nat$ and $ \alpha_1,\ldots,\alpha_{k_{l+1}-1}\in \Sigma$ such that 
$w_i =\alpha_{k_i}\alpha_{k_i+1}\ldots\alpha_{k_{i+1}-1}$ for all $i=1,\ldots,l$. 

If $\bw\in W^{\xi}(\Sigma)(s)$, then $s \propto w_1$ and $(s,w_1-s,\ldots, w_l)\in W^{\xi}(\Sigma)$. In case  $l>1$ we have that $n<k_2$ and that 
$\{n, k_2,\ldots, k_l \}\in \A_\xi$. So, $\{k_2,\ldots, k_l \}\in \A_{\xi_n}$ 
and consequently $(w_1,w_2,\ldots,w_l)\in W^{\xi_n}(\Sigma)\cap (W^{<\omega}(\Sigma) - s)$. In case 
$\bw=(w_1)\in W^{\xi}(\Sigma)(s)$, we have that $s \propto w_1$ and $(s,w_1-s)\in W^{\xi}(\Sigma)$. Thus $\{n\}\in \A_\xi$ and consequently $\emptyset\in \A_{\xi_n}$. This implies $\xi_n=0$ and indeed $\bw\in W^{\xi_n}(\Sigma)\cap (\W^{<\omega}(\Sigma) - s)$.  

If $\bw=(w_1,\ldots,w_l)\in W^{\xi_n}(\Sigma)\cap (W^{<\omega}(\Sigma) - s)$ and $l>1$, then $\{k_2,\ldots, k_l \} \in \A_{\xi_n} \cap {[\{n+1,n+2,\ldots \}]^{<\omega} \subseteq \A_\xi (n)}$. 
Hence, $\{n,k_2,\ldots, k_l\}\in \A_\xi$, so $(s,w_1-s,\ldots,w_k)\in W^{\xi}(\Sigma)$ and consequently $\bw\in W^{\xi}(\Sigma)(s)$. If $\bw=(w_1)\in W^{\xi_n}(\Sigma)\cap (W^{<\omega}(\Sigma) - s)$, then 
$\xi_n=0$ and consequently $\{n\}\in \A_\xi$. Hence $(s,w_1-s)\in W^{\xi}(\Sigma)$ and $\bw\in W^{\xi}(\Sigma)(s)$. 
\end{proof}

Now we can state and prove the main theorem of this section. 

\begin{thm}
[\textsf{A partition theorem on Schreier sets of words}]
\label{xi Carlson}
Let $\xi$ be a countable ordinal and $\Sigma$ a finite non-empty alphabet. For every $\G\subseteq W^{<\omega}(\Sigma)$, $\F\subseteq W^{<\omega}(\Sigma ; \upsilon)$ and every infinite sequence $\vec{w} \in W^\omega (\Sigma ; \upsilon)$ of variable words over $\Sigma$ there exists a variable reduction $\vec{u} \prec \vec{w} $ of $\vec{w}$ over $\Sigma$ such that :
\begin{itemize}
\item[{}] either $W^{\xi}(\Sigma)\cap RW^{<\omega}(\vec{u})\subseteq \G$ or 
$W^{\xi}(\Sigma)\cap RW^{<\omega}(\vec{u}) \subseteq W^{<\omega}(\Sigma)\setminus \G$, and
\end{itemize}
\begin{itemize}
\item[{}] either $ W^{\xi}(\Sigma ; \upsilon)\cap VRW^{<\omega}(\vec{u})\subseteq \F$ or 
$\W^{\xi}(\Sigma ; \upsilon)\cap VRW^{<\omega}(\vec{u})\subseteq \W^{<\omega}(\Sigma ; \upsilon)\setminus \F$. 
\end{itemize}
\end{thm}
For the proof of this partition theorem we will make use of a diagonal argument, contained in the following lemmas. 

\begin{lem}\label{lem:Carlson-Ramsey}
Let $\vec{w} = (w_n)_{n\in\nat} \in W^\omega (\Sigma ; \upsilon)$ be an infinite sequence of variable words over the alphabet $\Sigma$ and $\Pi_1 = \{(s,\vec{u}): s\in RW(\vec{w})$ and $\vec{u}\prec \vec{w} \setminus s \}$.\newline
If a subset $\R$ of $\Pi_1$ satisfies:
\begin{itemize}
\item[(i)] for every $(s,\vec{u})\in\Pi_1$ there exists $(s,\vec{u}_1)\in\R$ with 
$\vec{u}_1 \prec \vec{u}$; and
\item[(ii)] for every $(s,\vec{u}_1) \in\R$ and $\vec{u}_2 \prec \vec{u}_1$ we have $(s,\vec{u}_2)\in\R$,
\end{itemize}
then there exists $\vec{u} \prec \vec{w}$ such that 
$(s,\vec{s})\in \R$ for all $s\in RW(\vec{u})$ and $\vec{s} \prec \vec{u} \setminus s$.
\end{lem}
\begin{proof} 
Let $u_0 = w_1$ and $\vec{u}_0 =\vec{w}$. 
According to conditions (i) and (ii), there exists $\vec{u}_1 = (u^1_n)_{n\in\nat} \in 
W^\omega (\Sigma ; \upsilon)$ such that $\vec{u}_1 \prec \vec{w} \setminus u_0$ and 
$(u_0(\alpha),\vec{u}_1)\in \R$ for every $\alpha\in \Sigma$. 
Let $u_1 = u^1_1$. Then 
$(u_0, u_1)\in VRW^{<\omega}(\vec{w})$. 
We assume now that there have been constructed 
$\vec{u}_1,\ldots,\vec{u}_n \in W^\omega (\Sigma ; \upsilon)$ and $u_0,u_1,\ldots,u_n \in W(\Sigma ; \upsilon)$ 
such that $(u_0,u_1,\ldots, u_n)\in VRW^{<\omega}(\vec{w})$, $u_i\in VRW(\vec{u}_i)$, 
$\vec{u}_i \prec \vec{u}_{i-1} \setminus u_{i-1}$ 
for each $1\le i\le n$ and $(s,\vec{u}_i)\in\R$ for all $s\in RW((u_0,\ldots,u_{i-1}))$ and $1\le i\le n$.

We will construct $\vec{u}_{n+1}$ and $u_{n+1}$. 
Let $ RW ((u_0,\ldots, u_n))= \{ s_1,\ldots, s_k\}$ for some $k\in \nat$. 
Then $(s_i,\vec{u})\in\Pi_1$ for every $\vec{u}\prec \vec{u}_n \setminus u_n$ and $i=1,\ldots,k$.
According to condition (i), there exist 
$\vec{u}_{n+1}^1,\ldots, \vec{u}_{n+1}^k \in W^\omega (\Sigma ; \upsilon)$ such that 
$\vec{u}_{n+1}^k \prec \cdots \prec \vec{u}_{n+1}^1 \prec \vec{u}_n \setminus u_n$ and 
$(s_i, \vec{u}_{n+1}^i)\in \R$ for every $1\le i\le k$. 
Set $\vec{u}_{n+1} = \vec{u}_{n+1}^k$. If $\vec{u}_{n+1} = (u^{n+1}_i)_{i\in\nat}$, then set $u_{n+1} = u^{n+1}_1$.  
Of course $\vec{u}_{n+1}\prec \vec{u}_n \setminus u_n$, $u_{n+1} \in VRW(\vec{u}_{n+1})$, 
$(u_0,u_1,\ldots, u_{n+1})\in VRW^{<\omega}(\vec{w})$ and, according to condition (ii), 
$(s_i, \vec{u}_{n+1})\in \R$ for all $1\le i\le k$. 

Set $\vec{u} = (u_0,u_1,u_2,\ldots) \in W^\omega (\Sigma ; \upsilon)$. Then 
$\vec{u} \prec \vec{w}$, since $(u_0,u_1,\ldots, u_n)\in VRW^{<\omega}(\vec{w})$ 
for every $n\in\nat$. 
Let $s\in RW(\vec{u})$ and $\vec{s} \prec \vec{u} \setminus s$. 
Then there exists $n\in\nat$ such that $s\in RW((u_0,u_1,\ldots,u_n))$. 
Thus $(s,\vec{u}_{n+1})\in \R$ and, according to (ii), $(s,\vec{u} \setminus s)\in\R$, since 
$\vec{u} \setminus s \prec \vec{u}_{n+1}$. So $(s,\vec{s}) \in\R$, since $\vec{s}\prec\vec{u} \setminus s$.
\end{proof}

\begin{lem}\label{lem:Carlson-Ramsey 2}
Let $\vec{w} = (w_n)_{n\in\nat} \in W^\omega (\Sigma ; \upsilon)$ be an infinite sequence of variable words over the alphabet $\Sigma$ and $\Pi_2 = \{(t,\vec{u}): t\in VRW(\vec{w})$ and $\vec{u}\prec \vec{w}\setminus t \}$.\newline
If a subset $\R$ of $\Pi_2$ satisfies:
\begin{itemize}
\item[(i)] for every $(t,\vec{u})\in\Pi_2$ there exists $(t,\vec{u}_1)\in\R$ with 
$\vec{u}_1 \prec \vec{u}$; and
\item[(ii)] for every $(t,\vec{u}_1) \in\R$ and $\vec{u}_2 \prec \vec{u}_1$ we have $(t,\vec{u}_2)\in\R$,
\end{itemize}
then there exists $\vec{u} \prec \vec{w}$ such that 
$(t,\vec{t})\in \R$ for all $t\in VRW(\vec{u})$ and $\vec{t} \prec \vec{u}\setminus t$.
\end{lem}

\begin{proof} 
Let $u_0 = w_1$ and $\vec{u}_0 =\vec{w}$. 
According to condition (i), there exists $\vec{u}_1 = (u^1_n)_{n\in\nat} \in 
W^\omega (\Sigma ; \upsilon)$ 
such that $\vec{u}_1 \prec \vec{w} \setminus u_0$ and 
$(u_0,\vec{u}_1)\in \R$. 
Let $u_1 = u^1_1$. Then 
$(u_0, u_1)\in VRW^{<\omega}(\vec{w})$. 
The proof can be continued analogously to the proof of Lemma~\ref{lem:Carlson-Ramsey}.
\end{proof}

We are now ready to prove Theorem~\ref{xi Carlson}.

\begin{proof}[Proof of Theorem~\ref{xi Carlson}] 
Let $\G\subseteq W^{<\omega}(\Sigma)$, $\F \subseteq W^{<\omega}(\Sigma ; \upsilon)$ and 
$\vec{w}=(w_n)_{n\in \nat} \in W^\omega (\Sigma ; \upsilon)$. For $\xi=0$ the theorem is valid, according to Carlson's theorem, Theorem~\ref{thm:ICarlson}. Let $\xi>0$ be a countable ordinal. Assume that the theorem is valid for every $\zeta <\xi$. 

For every reduced word $s\in RW(\vec{w})$ of $\vec{w}$ over $\Sigma$ and every variable reduction 
$\vec{u} = (u_n)_{n\in \nat}\prec \vec{w}\setminus s $ of $\vec{w}\setminus s$ over $\Sigma$ is defined the variable reduction $\vec{u}_s = (s \ast u_1,u_2,\ldots) \prec \vec{w}$ of $\vec{w}$ over $\Sigma$. So, we can define the following set

 $\R_1= \{(s,\vec{u}): s\in RW(\vec{w}), \vec{u}\prec \vec{w}\setminus s $ and
\begin{itemize}
\item[{}] either $W^{\xi}(\Sigma)(s)\cap RW^{<\omega}(\vec{u}_s)\subseteq (\G\cap RW^{<\omega}(\vec{w}))(s)$ 
\item[{}] or $W^{\xi}(\Sigma)(s)\cap RW^{<\omega}(\vec{u}_s) \subseteq W^{<\omega}(\Sigma)\setminus (\G\cap RW^{<\omega}(\vec{w}))(s)\}.$
\end{itemize}
Of course $\R_1\subseteq \Pi_1 = \{(s,\vec{u}): s\in RW(\vec{w})$ and $\vec{u}\prec \vec{w} \setminus s \}$ and  obviously $\R_1$ satisfies the conditions (ii) of Lemma~\ref{lem:Carlson-Ramsey}. We will prove that $\R_1$ satisfies also the condition (i) of Lemma~\ref{lem:Carlson-Ramsey}. 

Let $(s,\vec{u})\in \Pi_1$. Then $s\in RW(\vec{w})\subseteq W(\Sigma)$, hence $s=\alpha_1\ldots\alpha_{n-1}$ for some $n\in \nat$ with $n>1$ and $\alpha_1,\ldots,\alpha_{n-1}\in \Sigma$. According to Proposition~\ref{justification 2}, there exists $\xi_n<\xi$ 
such that $W^{\xi}(\Sigma)(s) = W^{\xi_n}(\Sigma)\cap (W^{<\omega}(\Sigma) - s)$. 

If $\vec{u} = (u_n)_{n\in \nat}\prec \vec{w}\setminus s $, then $\vec{u}_s = (s \ast u_1,u_2,\ldots) \prec \vec{w}$. Using the induction hypothesis, there exists a variable reduction $\vec{u}^1 = (u^1_n)_{n\in \nat} \prec \vec{u}_s$ of $\vec{u}_s$ over $\Sigma$ such that  
\begin{itemize}
\item[{}] either $W^{\xi_n}(\Sigma)\cap RW^{<\omega}(\vec{u}^1)\subseteq (\G\cap RW^{<\omega}(\vec{w}))(s)$ 
\item[{}] or $\W^{\xi_n}(\Sigma)\cap R^{<\omega}(\vec{u}^1) \subseteq W^{<\omega}(\Sigma)\setminus (\G\cap RW^{<\omega}(\vec{w}))(s)$.
\end{itemize}
Then
\begin{itemize}
\item[{}] either $W^{\xi}(\Sigma)(s)\cap RW^{<\omega}(\vec{u}^1)\subseteq (\G\cap RW^{<\omega}(\vec{w}))(s)$ 
\item[{}] or $W^{\xi}(\Sigma)(s)\cap RW^{<\omega}(\vec{u}^1) \subseteq W^{<\omega}(\Sigma)\setminus (\G\cap RW^{<\omega}(\vec{w}))(s)$.
\end{itemize}
Since $\vec{u}^1 = (u^1_n)_{n\in \nat} \prec \vec{u}_s$, we have that $s \propto u^1_1$ so we set $\vec{u}_1 =  (u^1_1-s,u^1_2,\ldots)$. Then $\vec{u}_1 \prec \vec{u} \prec \vec{w}\setminus s$ and $(\vec{u}_1)_s = \vec{u}^1 $. Thus $(s,\vec{u}_1)\in\R_1$. Hence, $\R_1$ satisfies the condition (i) of Lemma~\ref{lem:Carlson-Ramsey}. 

According to Lemma~\ref{lem:Carlson-Ramsey}, there exists $\vec{w}_1  = (w^1_n)_{n\in \nat} \prec \vec{w}$ such that 
$(s,\vec{s})\in \R_1$ for all $s\in RW(\vec{w}_1)$ and $\vec{s} \prec \vec{w}_1\setminus s$. Thus, 
for every $s\in RW(\vec{w}_1)$ and $\vec{v} = (v_n)_{n\in \nat} \prec \vec{w}_1- s$, setting 
$\vec{v}_1 = (v_1-s,v_2,\ldots)$ we have that $(s,\vec{v}_1)\in\R_1$ and, since $(\vec{v}_1)_s = \vec{v}$ we have that
\begin{itemize}
\item[{}] either $W^{\xi}(\Sigma)(s)\cap RW^{<\omega}(\vec{v})\subseteq (\G\cap RW^{<\omega}(\vec{w}))(s)$ 
\item[{}] or $W^{\xi}(\Sigma)(s)\cap RW^{<\omega}(\vec{v}) \subseteq W^{<\omega}(\Sigma)\setminus (\G\cap RW^{<\omega}(\vec{w}))(s).$
\end{itemize}

Now, defining analogously for every variable reduced word $ t\in VRW(\vec{w}_1)$ of $\vec{w}_1$ over $\Sigma$ and every variable reduction 
$\vec{u} = (u_n)_{n\in \nat}\prec \vec{w}_1 \setminus t $ of $\vec{w}_1 \setminus t$ over $\Sigma$ the variable reduction 
$\vec{u}_t = (t \ast u_1,u_2,\ldots) \prec \vec{w}_1$ of $\vec{w}_1$ over $\Sigma$, we can define the set

$\R_2= \{(t,\vec{u}): t\in VRW(\vec{w}_1), \vec{u}\prec \vec{w}_1 \setminus t $ and
\begin{itemize}
\item[{}] either $ W^{\xi}(\Sigma ; \upsilon)(t)\cap VRW^{<\omega}(\vec{u}_t)\subseteq (\F\cap VRW^{<\omega}(\vec{w}))(t)$ 
\item[{}] or $W^{\xi}(\Sigma ; \upsilon)(t)\cap VRW^{<\omega}(\vec{u}_t)\subseteq W^{<\omega}(\Sigma ; \upsilon)\setminus (\F\cap VRW^{<\omega}(\vec{w}))(t)\}$. 
\end{itemize}
Then $\R_2\subseteq \Pi_2 = \{(t,\vec{u}): t\in VRW(\vec{w}_1)$ and $\vec{u}\prec \vec{w}_1\setminus t \}$ and $\R_2$ satisfies the condition (ii) of Lemma~\ref{lem:Carlson-Ramsey 2}.

Let $(t,\vec{u})\in \Pi_2$. Then $ t\in VRW(\vec{w}_1)$ and $t=\beta_1\ldots\beta_{n-1} \in W(\Sigma ; \upsilon)$ for some $n\in \nat$ with $n>1$ and $\beta_1,\ldots,\beta_{n-1} \in \Sigma \cup \{\upsilon\}$. According to Proposition~\ref{justification 2}, there exists $\xi_n<\xi$ 
such that $W^{\xi}(\Sigma ; \upsilon)(t) = W^{\xi_n}(\Sigma ; \upsilon)\cap (W^{<\omega}(\Sigma ; \upsilon) - t)$.

If $\vec{u} = (u_n)_{n\in \nat}\prec \vec{w}_1\setminus t $, then $\vec{u}_t = (t \ast u_1,u_2,\ldots) \prec \vec{w}_1$. Using the induction hypothesis, there exists a variable reduction $\vec{u}^1 = (u^1_n)_{n\in \nat} \prec \vec{u}_t$ of $\vec{u}_t$ over $\Sigma$ such that 
\begin{itemize}
\item[{}] either $ W^{\xi_n}(\Sigma ; \upsilon)\cap VRW^{<\omega}(\vec{u}^1)\subseteq (\F\cap VRW^{<\omega}(\vec{w}))(t)$ 
\item[{}] or $W^{\xi_n}(\Sigma ; \upsilon)\cap VRW^{<\omega}(\vec{u}^1)\subseteq W^{<\omega}(\Sigma ; \upsilon)\setminus (\F\cap VRW^{<\omega}(\vec{w}))(t)$. 
\end{itemize}
Then
\begin{itemize}
\item[{}] either $ W^{\xi}(\Sigma ; \upsilon)(t)\cap VRW^{<\omega}(\vec{u}^1)\subseteq (\F\cap VRW^{<\omega}(\vec{w}))(t)$ 
\item[{}] or $W^{\xi}(\Sigma ; \upsilon)(t)\cap VRW^{<\omega}(\vec{u}^1)\subseteq W^{<\omega}(\Sigma ; \upsilon)\setminus (\F\cap VRW^{<\omega}(\vec{w}))(t)$. 
\end{itemize}
Seting $\vec{u}_1 = (u^1_1-t,u^1_2,\ldots)$ we have that $\vec{u}_1 \prec \vec{u} \prec \vec{w}_1\setminus t$ and that 
$(t,\vec{u}_1)\in\R_2$. Hence, $\R_2$ satisfies also the condition (i) of Lemma~\ref{lem:Carlson-Ramsey 2} (replacing  $\vec{w}$ by $\vec{w}_1$).

According to Lemma~\ref{lem:Carlson-Ramsey 2}, there exists $\vec{w}_2  = (w^2_n)_{n\in \nat} \prec \vec{w}_1 \prec \vec{w}$ such that $(t,\vec{t})\in \R_2$ for all $t\in VRW(\vec{w}_2)$ and $\vec{t} \prec \vec{w}_2\setminus t$. Hence, for every $s\in RW(\vec{w}_2)\subseteq RW(\vec{w}_1)$, $t\in VRW(\vec{w}_2)$ and $\vec{v}_1 \prec \vec{w}_2- s\prec \vec{w}_1- s$, $\vec{v}_2 \prec \vec{w}_2- t$ we have
\begin{itemize}
\item[{}] either $W^{\xi}(\Sigma)(s)\cap RW^{<\omega}(\vec{v}_1)\subseteq (\G\cap RW^{<\omega}(\vec{w}))(s)$ 
\item[{}] or $W^{\xi}(\Sigma)(s)\cap RW^{<\omega}(\vec{v}_1) \subseteq W^{<\omega}(\Sigma)\setminus (\G\cap RW^{<\omega}(\vec{w}))(s);$ and
\end{itemize}
\begin{itemize}
\item[{}] either $ W^{\xi}(\Sigma ; \upsilon)(t)\cap VRW^{<\omega}(\vec{v}_2)\subseteq (\F\cap VRW^{<\omega}(\vec{w}))(t)$ 
\item[{}] or $W^{\xi}(\Sigma ; \upsilon)(t)\cap VRW^{<\omega}(\vec{v}_2)\subseteq W^{<\omega}(\Sigma ; \upsilon)\setminus (\F\cap VRW^{<\omega}(\vec{w}))(t)$. 
\end{itemize}
Let 
\begin{itemize}
\item[{}] $\G_1 =\{s\in RW(\vec{w}_2) : W^{\xi}(\Sigma)(s)\cap RW^{<\omega}(\vec{w}_2- s)\subseteq (\G\cap RW^{<\omega}(\vec{w}))(s)\}$ and 
\item[{}] $\F_1 = \{t\in VRW(\vec{w}_2) : W^{\xi}(\Sigma ; \upsilon)(t)\cap VRW^{<\omega}(\vec{w}_2- t)\subseteq (\F\cap VRW^{<\omega}(\vec{w}))(t)\}$. 
\end{itemize}
We use the induction hypothesis for $\xi=0$ (Theorem~\ref{thm:ICarlson}). Then, there exists a variable reduction $\vec{u} \prec \vec{w}_2$ of $\vec{w}_2$ such that : 
\begin{itemize}
\item[{}] either $RW(\vec{u}) \subseteq \G_1$ or 
$RW(\vec{u})\subseteq W(\Sigma )\setminus \G_1$; and
\end{itemize}
\begin{itemize}
\item[{}] either $VRW(\vec{u}) \subseteq \F_1$ or 
$VRW(\vec{u})\subseteq W(\Sigma ; \upsilon)\setminus \F_1$.
\end{itemize}
Since $\vec{u} \prec \vec{w}_2$, we have that $RW(\vec{u})\subseteq RW(\vec{w}_2)$ and 
$VRW(\vec{u})\subseteq VRW(\vec{w}_2)$. Thus
\begin{itemize}
\item[{}] either $W^{\xi}(\Sigma)(s)\cap RW^{<\omega}(\vec{u}- s)\subseteq (\G\cap RW^{<\omega}(\vec{w}))(s)$ for every 
$s\in RW(\vec{u})$
\item[{}] or $W^{\xi}(\Sigma)(s)\cap RW^{<\omega}(\vec{u}- s)\subseteq W^{<\omega}(\Sigma)\setminus (\G\cap RW^{<\omega}(\vec{w}))(s)$ for every $s\in RW(\vec{u})$; and
\end{itemize}
\begin{itemize}
\item[{}] either $W^{\xi}(\Sigma ; \upsilon)(t)\cap VRW^{<\omega}(\vec{u}- t)\subseteq (\F\cap VRW^{<\omega}(\vec{w}))(t)$ for every $t\in VRW(\vec{u})$ 
\item[{}] or $W^{\xi}(\Sigma ; \upsilon)(t)\cap VRW^{<\omega}(\vec{u}- t)\subseteq W^{<\omega}(\Sigma ; \upsilon)\setminus (\F\cap VRW^{<\omega}(\vec{w}))(t)$ for every $t\in VRW(\vec{u})$.
\end{itemize}
Hence,
\begin{itemize}
\item[{}] either $W^{\xi}(\Sigma)\cap RW^{<\omega}(\vec{u})\subseteq \G$ or 
$W^{\xi}(\Sigma)\cap RW^{<\omega}(\vec{u}) \subseteq W^{<\omega}(\Sigma)\setminus \G$, and
\end{itemize}
\begin{itemize}
\item[{}] either $ W^{\xi}(\Sigma ; \upsilon)\cap VRW^{<\omega}(\vec{u})\subseteq \F$ or 
$W^{\xi}(\Sigma ; \upsilon)\cap VRW^{<\omega}(\vec{u})\subseteq W^{<\omega}(\Sigma ; \upsilon)\setminus \F$. 
\end{itemize}
\end{proof}

We next give a more general statement of Theorem~\ref{xi Carlson}.

\begin{thm}
\label{xi-Furst-Kat}
Let $\xi$ be a countable ordinal, and $\vec{w}_0 \in \W^\omega (\Sigma ; \upsilon)$ be an infinite sequence of variable words over a finite, non-empty alphabet $\Sigma$. For any finite colorings $\chi_1 : RW^{\xi}(\vec{w}_0)\rightarrow \{1,\ldots,r_1\}$ and $\chi_2 : VRW^{\xi}(\vec{w}_0) \rightarrow \{1,\ldots,r_2 \}$ of the sets $RW^{\xi}(\vec{w}_0)$ and $VRW^{\xi}(\vec{w}_0)$ respectively and any variable reduction $\vec{w}\prec \vec{w}_0$ of $\vec{w}_0$ over $\Sigma$ there exists a variable reduction $\vec{u}\prec \vec{w}$ of $\vec{w}$ over $\Sigma$ such that all the finite reductions of $\vec{w}$ over $\Sigma$ in the set $RW^{\xi}(\vec{w}_0)$ are monochromatic under $\chi_1$ and all the finite variable reductions of $\vec{w}$ over $\Sigma$ in the set $VRW^{\xi}(\vec{w}_0)$ are monochromatic under $\chi_2$.
\end{thm}
\begin{proof}
Let $f : W^{<\omega}(\Sigma\cup\{\upsilon\}) \rightarrow RW^{<\omega}(\vec{w}_0)\cup VRW^{<\omega}(\vec{w}_0)$ with $f(\bs)=\vec{w}_0[\bs]$. Given the finite colorings $\chi_1 : RW^{\xi}(\vec{w}_0)\rightarrow \{1,\ldots,r_1\}$ and 
$\chi_2 : VRW^{\xi}(\vec{w}_0)\rightarrow \{1,\ldots,r_2\}$ are defined the finite colorings $\psi_1 : W^{<\omega}(\Sigma)\rightarrow \{1,\ldots,r_1\}$ where $\psi_1(\bs)= \chi_1(f(\bs))$ in case $\bs\in W^{\xi}(\Sigma)$ and $\psi_1(\bs)=1$ otherwise and $\psi_2 : W^{<\omega}(\Sigma ; \upsilon)\rightarrow \{1,\ldots,r_2\}$ where $\psi_2(\bt)= \chi_2(f(\bt))$ in case $\bt\in W^{\xi}(\Sigma ; \upsilon)$ and 
$\psi_2(\bt)=1$ otherwise. 

For a given $\vec{w}\prec \vec{w}_0$ there exists $\vec{t}\in \W^{\omega}(\Sigma ; \upsilon)$ such that $\vec{w}=\vec{w}_0[\vec{t}]$. According to Theorem~\ref{xi Carlson}, there exists a variable reduction $\vec{t}_1 \prec \vec{t} $ of $\vec{t}$ over $\Sigma$ such that the set $W^{\xi}(\Sigma)\cap RW^{<\omega}(\vec{t}_1)$ is 
monochromatic under $\psi_1$ and the set $ W^{\xi}(\Sigma ; \upsilon)\cap VRW^{<\omega}(\vec{t}_1)$ is 
monochromatic under $\psi_2$. Set $\vec{u}=\vec{w}_0[\vec{t}_1]\prec\vec{w}$. Then the set $RW^{\xi}(\vec{w}_0)\cap RW^{<\omega}(\vec{u})$ is monochromatic under $\chi_1$ and $ VRW^{\xi}(\vec{w}_0)\cap VRW^{<\omega}(\vec{u})$ is monochromatic under $\chi_2$. 
\end{proof}

We recall that in case $\vec{w}_0 = \vec{e}=(e_n)_{n\in\nat}$ with $e_n=\upsilon$ for every $n\in\nat$ all the  infinite sequences of variable words over $\Sigma$ are variable reductions of $\vec{w}_0$ over $\Sigma$ and that 
$RW^{\xi}(\vec{w}_0)=W^{\xi}(\Sigma)$, $VRW^{\xi}(\vec{w}_0)= W^{\xi}(\Sigma ; \upsilon)$ for every $0\leq \xi<\omega_1$. In this case Theorem~\ref{xi-Furst-Kat} is indentified with Theorem A referred to the introduction.

For $k\in\nat$ and $\vec{u} \in \W^\omega (\Sigma ; \upsilon)$ we have that 
$W^{k}(\Sigma)\cap RW^{<\omega}(\vec{u})= RW^{k}(\vec{u})$ and $W^{k}(\Sigma ; \upsilon)\cap VRW^{<\omega}(\vec{u})= VRW^{k}(\vec{u})$, hence Theorem~\ref{xi-Furst-Kat} in case $\xi=k\in\nat$ implies Theorem~\ref{IFurst-Kat} which essencially has be proved by Furstenberg and Katznelson in \cite{FK1} (Theorems 2.7 and 3.1).

The following theorem is a finitary consequence of Theorem~\ref{xi-Furst-Kat}. It follows from Theorem~\ref{xi-Furst-Kat} using a compactness argument. We will need the following notation to state it. For a word $w=\alpha_1\ldots\alpha_l$ over an alphabet $\Sigma$ let $l$ be the length of $w$. We denote by 
$W_l(\Sigma)$ the set of all words over $\Sigma$ with length $l$. For a countable ordinal $\xi$, we denote by  $W^{\xi}_M (\Sigma)$ the set of all finite sequences of words in $W^{\xi}(\Sigma)$ such that the sum of the lengths of their words is equal to $M$.

\begin{thm}
[\textsf{Extended Hales-Jewett theorem}]
\label{xi-Ha-Je}
For every $r, n, k \in\nat$, $\Sigma$ a finite, non-empty alphabet of $k$ elements and $\xi$ a  countable ordinal there exists $M=M( r,n,k,\xi )$ such that for every $r$-coloring of $W^{\xi}_M (\Sigma)$ there exists a finite sequence $\bw=(w_1,\ldots,w_n)$ of variable words over $\Sigma$ all of whose the finite reductions over $\Sigma$ in  $W^{\xi}_M (\Sigma)$ are monochromatic.
\end{thm}

The classical Hales-Jewett theorem (\cite{HJ}), is a trivial consequence of the case $\xi = 0$, $n = 1$. Since van der Waerden's theorem (\cite{vdW}) may be obtained as a corollary of the Hales-Jewett theorem, Theorem~\ref{xi-Ha-Je} can be used to obtain a corresponding extention of van der Waerden's theorem. 

Furstenberg and Katznelson in \cite{FK1} introduced the notion of a $k$-dimensional combinatorial subspace of $W(\Sigma)$ for $k\in\nat$ and proved (in Theorem 3.1) a partition theorem about these combinatorial subspaces.  Theorem~\ref{xi Carlson} implies an extension of this partition theorem to every countable ordinal. Let give the neccesary notation.

Let $\Sigma$ be a finite, non-empty alphabet. A \textit{finite-dimensional combinatorial subspace} $[\bw]$ of $W(\Sigma)$ is defined by a finite sequence $\bw=(w_1,\ldots,w_k)\in \W^{<\omega}(\Sigma; \upsilon)$ of variable words over $\Sigma$ as follows:
\begin{center}
$[\bw]=RW((w_1,\ldots,w_k)) = \{w_1(\alpha_1)\ast \ldots \ast w_k(\alpha_k): \alpha_1,\ldots,\alpha_k \in \Sigma \}$.
\end{center}

In the same way an \textit{infinite-dimensional combinatorial subspace} $[\vec{w}]$ of $W(\Sigma)$ is defined by an infinite sequence $\vec{w}=(w_n)_{n\in\nat} \in \W^{\omega}(\Sigma; \upsilon)\}$ as follows:
\begin{center}
$[\vec{w}]=RW(\vec{w})=\{w_1(\alpha_1)\ast \ldots \ast w_k(\alpha_k): k\in \nat, \alpha_1,\ldots,\alpha_k \in \Sigma\}$.
\end{center}

A finite (or infinite)-dimensional combinatorial subspace of $W(\Sigma)$ contained in an infinite-dimensional combinatorial subspace $[\vec{w}]$ of $W(\Sigma)$ is called \textit{a finite (or infinite)-dimensional combinatorial subspace} of $[\vec{w}]$. It is not hard to check that a finite-dimensional combinatorial subspaces of $[\vec{w}]$ is of the form  $[\bu]$, where $\bu \in VRW^{<\omega}(\vec{w})$ and that an infinite-dimensional combinatorial subspaces of $[\vec{w}]$ is of the form  $[\vec{u}]$, where $\vec{u} \in VRW^{\omega}(\vec{w})$. 

\begin{defn} 
\label{dimcomb}
Let $\xi$ be a countable ordinal. A \textit{$\xi$-combinatorial subspace} $[\bw]$ of $W(\Sigma)$ is a finite-dimensional combinatorial subspace of $W(\Sigma)$ such that $\bw \in W^{\xi}(\Sigma; \upsilon)$ and a  \textit{$\xi$-combinatorial subspace} $[\bw]$ of an infinite-dimensional combinatorial subspace $[\vec{w}]$ of $W(\Sigma)$ is a finite-dimensional combinatorial subspace of $[\vec{w}]$ such that $\bw \in VRW^{\xi}(\vec{w})$. 
\end{defn}

The class of $k$-combinatorial subspaces of $W(\Sigma)$, for $k\in \nat$, coincites with the class of $k+1$-dimensional combinatorial subspaces of $W(\Sigma)$, while the class of $\xi$-combinatorial subspaces of $W(\Sigma)$, for 
a countable ordinal $\xi\geq \omega$, contains finite-dimensional combinatorial subspaces of  $W(\Sigma)$ of arbitrary large finite dimentions. Also, observe that although for $k\in \nat$ the $k$-combinatorial subspaces of $[\vec{w}]$ are exactly the $k$-combinatorial subspaces of $W(\Sigma)$ contained in $[\vec{w}]$, for a countable ordinal $\xi\geq \omega$ it is not always true that every $\xi$-combinatorial subspaces of $[\vec{w}]$ is a $\xi$-combinatorial subspaces of $W(\Sigma)$, since it is not true that $d_{\vec{w}}(\bw)= d(\bw)$ for every $\bw\in R^{<\omega}(\vec{w})$.  

We will state now a corollary of Theorem~\ref{xi Carlson} which extents Theorem 3.1 in \cite{FK2},  corresponding to finite ordinals $\xi < \omega$, to every countable ordinal $\xi$. 

\begin{cor}
[\textsf{Combinatorial subspaces partition theorem}]
\label{Com Furst-Kat}
Let $\xi$ be a countable ordinal. For any finite coloring of the set $CS^{\xi}(\Sigma)$ of all $\xi$-combinatorial subspaces of $W(\Sigma)$ and any infinite-dimensional combinatorial subspace $[\vec{w}]$ of $W(\Sigma)$, there exists an infinite-dimensional combinatorial subspace $[\vec{u}]$ of $[\vec{w}]$ such that all the $\xi$-combinatorial subspaces of $W(\Sigma)$ contained in $[\vec{u}]$ are monochromatic.
\end{cor} 
\begin{proof}
Given the finite coloring $\chi : CS^{\xi}(\Sigma)\rightarrow \{1,\ldots,r\}$, is defined the finite coloring 
$\psi : \W^{\xi}(\Sigma; \upsilon)\rightarrow \{1,\ldots,r\}$ with $\psi(\bs)= \chi([\bs])$. Apply Theorem~\ref{xi-Furst-Kat} for $\vec{w}_0 =(e_n)_{n\in\nat}$ with $e_n=\upsilon$ for every $n\in \nat$. Then 
for any $\vec{w}\in \W^\omega (\Sigma ; \upsilon)$, there exists a variable reduction $\vec{u}$ of $\vec{w}$ over $\Sigma$ such that all the elements of the set $W^{\xi}(\Sigma ; \upsilon)\cap VRW^{<\omega}(\vec{u})$ are $\psi$-monochromatic. Hence, for any infinite-dimensional combinatorial subspace $[\vec{w}]$ of $W(\Sigma)$, there exists an infinite-dimensional combinatorial subspace $[\vec{u}]$ of $[\vec{w}]$ such that all the $\xi$-combinatorial subspaces of $W(\Sigma)$ contained in $[\vec{u}]$ are $\chi$-monochromatic. 
\end{proof} 

\begin{cor}
\label{Furst-Kat2}
Let $\xi$ be a countable ordinal, and $\vec{w}_0 \in \W^\omega (\Sigma ; \upsilon)$. For any finite coloring of all $\xi$-combinatorial subspaces of $[\vec{w}_0]$ and any infinite-dimensional combinatorial subspace $[\vec{w}]$ of $\vec{w}_0$, there exists an infinite-dimensional combinatorial subspace $[\vec{u}]$ of $[\vec{w}]$ such that all the $\xi$-combinatorial subspaces of $[\vec{w}_0]$ contained in $[\vec{u}]$ are monochromatic.
\end{cor} 

Using the previous terminology we obtain a generalization of Hales-Jewett theorem to higher dimensions, as a consequence of Theorem~\ref{xi-Ha-Je}.
\begin{cor}
\label{xi-Ha-Je1}
For every $r, n, k \in\nat$, $\Sigma$ a finite, non-empty alphabet of $k$ elements and $\xi$ a  countable ordinal there exists $M=M( r,n,k,\xi )$ such that for any $r$-coloring of $W_M (\Sigma)$ there exists a monochromatic $\xi$-combinatorial subspace of $W(\Sigma)$.
\end{cor}

\section{Basic properties of the Schreier type families of the finite sequences of words}

This section is preparatory for the results of sections 4 and 5. We prove here (a) the thiness of the Schreier-type families of words $W^{\xi}(\Sigma)$ and variable words $W^{\xi}(\Sigma; \upsilon)$ (Proposition~\ref{prop:thinfamily}), and (b) the canonical representation of every (infinite or finite) sequence of (variable) words over $\Sigma$ with respect to the Schreier-type families  (Proposition~\ref{prop:canonicalrep}). Furthermore we introduce the (strong) Cantor-Bendixson index of a hereditary subfamily of the family of the finite sequences of (variable) words  (Definition~\ref{def:Cantor-Bendix}), and we prove that the index of the hereditary family generated by the $\xi$-Schreier-type family of finite sequences of words is $\xi+1$ for every countable ordinal $\xi$ (Proposition~\ref{prop:Cantor-Bendix}). In addition, in Theorem~\ref{thm:treeRam}, we strengthen Theorem~\ref{xi Carlson} in case the partition family is (not an arbitrary family but) a tree.

\begin{defn}\label{def:Fthin1}
Let $\Sigma$ be a finite, non empty alphabet and $\F\subseteq W^{<\omega}(\Sigma \cup \{\upsilon\})$ be a family of finite  sequences of words over $\Sigma \cup \{\upsilon\}$.
\begin{itemize}
\item[(i)] $\F$ is {\em thin\/} if there are no elements $\bs=(s_1,\ldots,s_k),\bt=(t_1,\ldots,t_k)\in\F$
with $\bs \propto \bt$ (which means that $k<l$ and $s_i=t_i$ for all $i=1,\ldots,k$).
\item[(ii)] $\F^* = \F \cup \{\bt \in W^{<\omega}(\Sigma \cup \{\upsilon\}): \bt\propto \bs$ for some 
$\bs\in \F\}\cup \{\emptyset\}$.
\item[(iii)] $\F$ is a {\em tree\/} if $\F^* = \F$.
\end{itemize}
\end{defn}

\begin{prop}\label{prop:thinfamily}
Let $\vec{w}=(w_n)_{n\in\nat} \in W^\omega (\Sigma ; \upsilon)$ be an infinite sequence of variable words over an alphabet $\Sigma$. The families $W^{\xi}(\Sigma ; \upsilon)$, $W^{\xi}(\Sigma)$, $VRW^{\xi}(\vec{w})$, $RW^{\xi}(\vec{w})$ are thin for every $\xi<\omega_1$.
\end{prop} 
\begin{proof} 
It follows from the fact that the families $\A_\xi$ of Schreier finite subsets of $\nat$ are thin (which means that if 
$s,t\in \A_\xi$ and $s$ is an initial segment of $t$, then $s=t$).
\end{proof}

\begin{prop}\label{prop:canonicalrep}
Let $\xi> 0$ be a countable ordinal number.

{\rm (i)} Every infinite sequence $\vec{s} = (s_n)_{n\in\nat}\in W^{\omega}(\Sigma \cup \{\upsilon\})$ of words over  $\Sigma\cup \{\upsilon\}$ has canonical representation with 
respect to $W^{\xi}(\Sigma\cup \{\upsilon\})$, which means that there exists a unique strictly increasing sequence $(m_n)_{n\in\nat}$ in $\nat$ such that $(s_1,\ldots,s_{m_1}) \in W^{\xi}(\Sigma\cup \{\upsilon\})$ and $(s^n, s_{m_{n-1}+1},\ldots, s_{m_n}) \in W^{\xi}(\Sigma\cup \{\upsilon\})$ for every $n > 1$, where $s^n =s_1\ast\ldots\ast s_{m_{n-1}} $. 

{\rm (ii)} Every nonempty finite sequence $\bs = (s_1,\ldots,s_k)\in W^{<\omega}(\Sigma \cup \{\upsilon\})$  
of words over  $\Sigma\cup \{\upsilon\}$ has canonical representation 
with respect to $W^{\xi}(\Sigma\cup \{\upsilon\})$, which means that either $\bs\in \big(W^{\xi}(\Sigma\cup \{\upsilon\})\big)^* \setminus W^{\xi}(\Sigma\cup \{\upsilon\})$ or there exist 
unique $n\in\nat$, and  $m_1, \ldots,m_n \in\nat$
with $m_1 < \ldots < m_n\leq k$ such that $(s_1,\ldots,s_{m_1})\in W^{\xi}(\Sigma\cup \{\upsilon\})$, $(s^n, s_{m_{n-1}+1},\ldots,s_{m_n}) \in W^{\xi}(\Sigma\cup \{\upsilon\})$ for every $n > 1$, where 
$s^n =s_1\ast\ldots\ast s_{m_{n-1}} $  and in case $m_n < k$,  
$(s^{n+1}, s_{{m_n}+1},\ldots,s_k)\in \big(W^{\xi}(\Sigma\cup \{\upsilon\})\big)^* \setminus W^{\xi}(\Sigma\cup \{\upsilon\})$ where $s^{n+1} = s_1\ast\ldots\ast s_{m_{n}} $.
\end{prop}
\begin{proof}
(i) Let 
$\xi>0$ and $\vec{s} = (s_n)_{n\in\nat}\in W^{\omega}(\Sigma \cup \{\upsilon\})$. Then the sequence 
$d((s_n)_{n\in\nat}) = (k_n)_{n\geq 2}$ of natural numbers has canonical representation with respect 
to $\A_\xi$, which means that there exists a unique strictly increasing sequence $(m_n)_{n\in\nat}$ in $\nat$ so that $(k_2,\ldots,k_{m_1})\in \A_\xi$ and $(k_{m_{n-1}+1},\ldots,k_{m_n})\in \A_\xi$ for every $n > 1$. Hence, 
$(s_1,\ldots,s_{m_1}) \in W^{\xi}(\Sigma\cup \{\upsilon\})$ and $(s^n, s_{m_{n-1}+1},\ldots, s_{m_n}) \in  W^{\xi}(\Sigma\cup \{\upsilon\})$ for every $n > 1$, where $s^n =s_1\ast\ldots\ast s_{m_{n-1}}$.

(ii) Let $\bs = (s_1,\ldots,s_k)\in W^{<\omega}(\Sigma \cup \{\upsilon\})$. Set $s_n = \upsilon$ for every $n\in \nat$ with $n>k$. The sequence $\vec{s} = (s_n)_{n\in\nat}\in W^{\omega}(\Sigma \cup \{\upsilon\})$ has canonical representation with respect to $W^{\xi}(\Sigma\cup \{\upsilon\})$, according to (i).
\end{proof}

According to Proposition~\ref{prop:canonicalrep}, every finite or infinite reduction (or variable reduction) of a sequence $\vec{w} = (w_n)_{n\in\nat}\in W^{\omega}(\Sigma ; \upsilon)$ has canonical representation with respect to $RW^{\xi}(\vec{w})$ (or to $VRW^{\xi}(\vec{w})$), for every $1\leq \xi<\omega_1$. For example $\b{u}=\vec{w}[\bs]\in RW^{<\omega}(\vec{w})$ has canonical representation with respect to $RW^{\xi}(\vec{w})$, as $\bs$ has canonical representation with respect to $W^{\xi}(\Sigma)$.

Now exploiting the canonical representation of every sequence of words over $\Sigma \cup \{\upsilon\}$ with 
respect to $W^{\xi}(\Sigma\cup \{\upsilon\})$ we will give alternative descriptions of the dichotomies described in Theorem~\ref{xi Carlson}.

\begin{prop}\label{prop:canonicalrep2}
Let $\G\subseteq W^{<\omega}(\Sigma)$, $\F \subseteq W^{<\omega}(\Sigma ; \upsilon)$ and let $\xi$ be a countable ordinal. Then, for every infinite sequence $\vec{u}=(u_n)_{n\in\nat}\in W^\omega (\Sigma ; \upsilon)$ of variable words over $\Sigma$ the following are equivalent:

{\rm (i)} $W^{\xi}(\Sigma)\cap RW^{<\omega}(\vec{u})\subseteq \G$ ( resp. $ W^{\xi}(\Sigma ; \upsilon)\cap VRW^{<\omega}(\vec{u})\subseteq \F$).

{\rm (ii)} For every variable reduction $\vec{u}_1 $ of $\vec{u}$ the unique initial segment $\bs = (u^1_1,\ldots,u^1_m)$  of $\vec{u}_1$ which is an element of $W^{\xi}(\Sigma ; \upsilon)$ satisfies the property  $(u^1_1(\alpha_1),\ldots,u^1_m(\alpha_m))\in \G$ for every $\alpha_1,\ldots,\alpha_m\in \Sigma$ (resp. the property $\bs\in \F$).

{\rm (iii)} Given any sequence $(\vec{u}_n)_{n\in\nat}$ of infinite sequences of variable words over $\Sigma$ such that $\vec{u}_1\prec\vec{u}$ and $\vec{u}_{n+1}\prec\vec{u}_{n}$ for every $n\in\nat$ and any 
$t_n \in VRW(\vec{u}_n)$ with $t_n \propto t_{n+1}$ for every $n\in\nat$, 
there exists $m\in\nat$ such that $(t^1_1(\alpha_1),\ldots,t^1_{m}(\alpha_{m})) \in W^{\xi}(\Sigma)\cap \G$ for every $\alpha_1,\ldots,\alpha_{m}\in\Sigma$, where $t^1_1=t_1$ and $t^1_i=t_i-t_{i-1}$ for $i=2,\ldots,m$ (resp. such that $(t^1_1,\ldots,t^1_{m}) \in W^{\xi}(\Sigma ; \upsilon)\cap \F$).
\end{prop}
\begin{proof} 
 (i) $\Rightarrow$ (ii). Let $\vec{u}_1 =(u^1_n)_{n\in\nat}$ be a variable reduction of $\vec{u}$. 
Using the canonical representation of $\vec{u}_1$ with respect to $W^{\xi}(\Sigma ; \upsilon)$  (Proposition~\ref{prop:canonicalrep}), there exists a unique initial segment $\bs=(u^1_1,u^1_2,\ldots,u^1_m)$ of $\vec{u}_1$ which is an element of $W^{\xi}(\Sigma ; \upsilon)$. According to (i), $(u^1_1(\alpha_1),\ldots,u^1_m(\alpha_m)))\in W^{\xi}(\Sigma)\cap RW^{<\omega}(\vec{u})\subseteq \G$ for every $\alpha_1,\ldots,\alpha_m\in \Sigma$ (resp. $\bs \in \F$).

(ii) $\Rightarrow$ (i). 
Let $\bs = (s_1,\ldots,s_m)\in W^{\xi}(\Sigma)\cap RW^{<\omega}(\vec{u})$. There exist unique sequences $1=k_1<\cdots<k_m<k_{m+1}\in \nat$ and $ \alpha_1,\ldots, \alpha_{k_{m+1}-1}\in \Sigma$ such that 
$s_i = u_{k_i}(\alpha_{k_i})\ast \ldots \ast u_{k_{i+1}-1}(\alpha_{k_{i+1}-1})$ for all 
$i=1,\ldots,m$. Set $u^1_i = u_{k_i}(\upsilon)\ast u_{k_i+1}(\alpha_{k_i+1})\ast \ldots \ast u_{k_{i+1}-1}(\alpha_{k_{i+1}-1})$ for all $i=1,\ldots,m$ and $u^1_{m+i}= u_{k_{m+1}-1+i}$ for every $i\in\nat$. Then 
the sequence $\vec{u}_1 = (u^1_n)_{n\in\nat}\in W^\omega (\Sigma ; \upsilon)$ is a variable reduction of $\vec{u}$ 
and $(u^1_1,\ldots,u^1_m)\in W^{\xi}(\Sigma ; \upsilon)$. According to (ii), we have that $\bs\in \G$. 

If $\bt = (t_1,\ldots,t_m)\in W^{\xi}(\Sigma ; \upsilon)\cap VRW^{<\omega}(\vec{u})$, then $\bt$ is the unique initial segment of a variable reduction $\vec{u}_1$ of $\vec{u}$, hence, according to (ii), $\bt\in \F$.     

(ii) $\Rightarrow$ (iii). Let a sequence $(\vec{u}_n)_{n\in\nat}$ of infinite sequences of variable words over $\Sigma$ such that $\vec{u}_1\prec\vec{u}$ and $\vec{u}_{n+1}\prec\vec{u}_{n}$ for every $n\in\nat$ and 
$t_n \in VRW(\vec{u}_n)$ with $t_n \propto t_{n+1}$ for every $n\in\nat$. The sequence 
$\vec{t}=(t^1_n)_{n\in\nat}$ with $t^1_1=t_1$ and $t^1_n=t_n-t_{n-1}$ for $n>1$ is a variable reduction of $\vec{u}$, hence, according to (ii), there exists $m\in\nat$ such that $(t^1_1,\ldots,t^1_m)\in W^{\xi}(\Sigma ; \upsilon)$ and 
$(t^1_1(\alpha_1),\ldots,t^1_{m}(\alpha_{m})) \in W^{\xi}(\Sigma)\cap \G$ for every $\alpha_1,\ldots,\alpha_{m}\in\Sigma$ (resp. and $(t^1_1,\ldots,t^1_{m}) \in \F$). 

(iii) $\Rightarrow$ (ii). Let a variable reduction $\vec{u}_1= (u^1_n)_{n\in\nat} $ of $\vec{u}$. Use (iii), seting $\vec{u}_n = \vec{u}_1$ for every $n\in\nat$ and $t_n=u^1_1\ast \ldots \ast u^1_n$ for 
every $n\in\nat$.
\end{proof}

We will give now an alternative description for the second horn of the dichotomy proved in Theorem~\ref{xi Carlson}, in case the partition family is a tree.

\begin{prop}\label{prop:tree}
Let $\G\subseteq W^{<\omega}(\Sigma)$ be a tree, $\F\subseteq W^{<\omega}(\Sigma ; \upsilon)$ be a tree and let $\xi$ be a countable ordinal. Then 
\newline
$W^{\xi}(\Sigma)\cap RW^{<\omega}(\vec{u}) \subseteq W^{<\omega}(\Sigma)\setminus \G$  
if and only if 
$\G\cap RW^{<\omega}(\vec{u}) \subseteq \big(W^{\xi}(\Sigma)\big)^* \setminus W^{\xi}(\Sigma)$ and 
\newline
$W^{\xi}(\Sigma ; \upsilon)\cap VRW^{<\omega}(\vec{u})\subseteq W^{<\omega}(\Sigma ; \upsilon)\setminus \F$
if and only if 
$\F\cap VRW^{<\omega}(\vec{u}) \subseteq \big(W^{\xi}(\Sigma ; \upsilon)\big)^* \setminus W^{\xi}(\Sigma ; \upsilon)$.
\end{prop}
\begin{proof} 
Let $W^{\xi}(\Sigma)\cap RW^{<\omega}(\vec{u}) \subseteq W^{<\omega}(\Sigma)\setminus \G$ 
and $\bs = (s_1,\ldots,s_l)\in \G\cap RW^{<\omega}(\vec{u})$. Since $\bs$ has canonical representation with respect to $W^{\xi}(\Sigma)$ (Proposition~\ref{prop:canonicalrep}), either $\bs\in \big(W^{\xi}(\Sigma)\big)^* \setminus W^{\xi}(\Sigma)$, as required, or there 
exists $\bs_1\in W^{\xi}(\Sigma)$ such that either $\bs_1 = \bs$ or $\bs_1 \propto \bs$. 
The second case is impossible, since then $\bs_1 \in \G\cap RW^{<\omega}(\vec{u}) \cap W^{\xi}(\Sigma)$; a contradiction to our assumption. Hence, $\G\cap RW^{<\omega}(\vec{u}) \subseteq \big(W^{\xi}(\Sigma)\big)^* \setminus W^{\xi}(\Sigma)$. Obviously,  
$W^{\xi}(\Sigma)\cap RW^{<\omega}(\vec{u}) \subseteq W^{<\omega}(\Sigma)\setminus \G$ 
if $\G\cap RW^{<\omega}(\vec{u}) \subseteq \big(W^{\xi}(\Sigma)\big)^* \setminus W^{\xi}(\Sigma)$. 

Analogously, can be proved that $W^{\xi}(\Sigma ; \upsilon)\cap VRW^{<\omega}(\vec{u})\subseteq W^{<\omega}(\Sigma ; \upsilon)\setminus \F$ 
if and only if 
$\F\cap VRW^{<\omega}(\vec{u}) \subseteq \big(W^{\xi}(\Sigma ; \upsilon)\big)^* \setminus W^{\xi}(\Sigma ; \upsilon)$.
\end{proof}

A consequence of Proposition~\ref{prop:tree} is the following stronger form of Theorem~\ref{xi Carlson} in case the partition families are trees.

\begin{thm}\label{thm:treeRam}
Let $\xi$ be a countable ordinal, $\Sigma$ a finite non-empty alphabet and $\G\subseteq W^{<\omega}(\Sigma)$, $\F\subseteq W^{<\omega}(\Sigma ; \upsilon)$ be trees. Then for every infinite sequence $\vec{w} \in W^\omega (\Sigma ; \upsilon)$ of variable words over $\Sigma$ there exists a variable reduction $\vec{u} \prec \vec{w} $ of $\vec{w}$ over $\Sigma$ such that :
\begin{itemize}
\item[{}] either $W^{\xi}(\Sigma)\cap RW^{<\omega}(\vec{u})\subseteq \G$ or 
$\G\cap RW^{<\omega}(\vec{u}) \subseteq \big(W^{\xi}(\Sigma)\big)^* \setminus W^{\xi}(\Sigma)$, and
\end{itemize}
\begin{itemize}
\item[{}] either $ W^{\xi}(\Sigma ; \upsilon)\cap VRW^{<\omega}(\vec{u})\subseteq \F$ or 
$\F\cap VRW^{<\omega}(\vec{u}) \subseteq \big(W^{\xi}(\Sigma ; \upsilon)\big)^* \setminus W^{\xi}(\Sigma ; \upsilon)$.
\end{itemize}
\end{thm}

We will define now a topology on the sets $W^{<\omega}(\Sigma)$, $W^{<\omega}(\Sigma ; \upsilon)$, $W^{\omega}(\Sigma)$, $W^{\omega}(\Sigma ; \upsilon)$. We set $D=\{(n,\alpha): n\in \nat, \alpha \in \Sigma\cup\{\upsilon\} \}$, note that $D$ is a countable set, and denote by $[D]^{<\omega}$ the set of all finite subsets of $D$. 

Each finite sequence $\bw$ of words over $\Sigma\cup\{\upsilon\}$ corresponds to a unique finite subset $\sigma(\bw)$ of $[D]^{<\omega}$ defined as follows: We set $\sigma(\emptyset)=\emptyset$. For $\bw=(w_1,\ldots,w_m)\in W^{<\omega}(\Sigma\cup\{\upsilon\})$ there exist $1=k_1<\cdots<k_m<k_{m+1}\in\nat$ and $ \alpha_i \in \Sigma\cup\{\upsilon\}$ for all $i=1,\ldots,k_{m+1}-1$ such that $w_i =\alpha_{k_i}\ldots\alpha_{k_{i+1}-1}$ for all $i=1,\ldots,m$, hence we set
\begin{center}
$\sigma(\bw)=\{ \{(k_j,\alpha_{k_j}) : j\in \nat, k_i\leq j\leq k_{i+1}-1 \} :  i\in\nat, i\leq m \}$.
\end{center}
Analogously, for $\vec{w}=(w_n)_{n\in\nat}\in W^{\omega}(\Sigma\cup\{\upsilon\})$ there exist $1=k_1<k_2<k_3<\cdots\in \nat$ and $ \alpha_i \in \Sigma\cup\{\upsilon\}$ for every $i\in\nat$ such that $w_n=\alpha_{k_n}\ldots\alpha_{k_{n+1}-1}$ for all $n\in \nat$, hence we set 
$\sigma(\vec{w})=\{ \{(k_j,\alpha_{k_j}) : j\in \nat, k_n\leq j\leq k_{n+1}-1 \} :  n\in\nat \}$.

We identify every sequence (finite or infinite) of words over $\Sigma\cup\{\upsilon\}$ with its characteristic function in $\{0,1\}^{[D]^{<\omega}}$ via the function: 
\begin{center}
$I : W^{<\omega}(\Sigma\cup\{\upsilon\})\cup W^{\omega}(\Sigma\cup\{\upsilon\}) \rightarrow \{0,1\}^{[D]^{<\omega}}$ 
\end{center}
\begin{center}
with $I(\bw)=x_{\sigma(\bw)} $ for $\bw\in W^{<\omega}(\Sigma\cup\{\upsilon\})$ and 
$I(\vec{w})=x_{\sigma(\vec{w})} $ for $\vec{w}\in W^{\omega}(\Sigma\cup\{\upsilon\})$,
\end{center}
Thus, identifying every finite sequence $\bw\in W^{<\omega}(\Sigma\cup\{\upsilon\})$ and every infinite sequence $\vec{w}\in W^{\omega}(\Sigma\cup\{\upsilon\})$ of words over $\Sigma\cup\{\upsilon\}$ with its characteristic function $x_{\sigma(\bw)} \in \{0,1\}^{[D]^{<\omega}}$ and $x_{\sigma(\vec{w})} \in \{0,1\}^{[D]^{<\omega}}$ respectively, we topologize the sets $W^{<\omega}(\Sigma)$, $W^{<\omega}(\Sigma ; \upsilon)$, $W^{\omega}(\Sigma)$, $W^{\omega}(\Sigma ; \upsilon)$ by the topology of pointwise convergence (equivalently 
by the relative product topology of $\{0,1\}^{[D]^{<\omega}}$). For example we say that a family $\F\subseteq W^{<\omega}(\Sigma ; \upsilon)$ is pointwise closed iff the family $\{x_{\sigma(\bw)} :\bw \in \F\}$ is closed in the 
topology of pointwise convergence or a family $\U\subseteq W^{\omega}(\Sigma)$ is pointwise closed iff 
$\{x_{\sigma(\vec{w})} :\vec{w}\in \U\}$ is pointwise closed in $\{0,1\}^{[D]^{<\omega}}$.

We next turn our attention to hereditary families of finite sequences of words. 
\begin{defn}\label{def:Hthin}
Let $\Sigma$ be a finite, non empty alphabet and $\F\subseteq W^{<\omega}(\Sigma ; \upsilon)$.
\begin{itemize}
\item[(i)] $\F_* = \{\bt \in W^{<\omega}(\Sigma ; \upsilon): \bt\in VRW^{<\omega}(\bs)$ 
for some $\bs \in \F^*\setminus \{\emptyset\} \}\cup\{\emptyset\}$.
\item[(ii)] $\F$ is {\em hereditary\/} if $\F_* = \F$.
\end{itemize}
\end{defn}

\begin{defn}\label{def:Herthin}
Let $\Sigma$ be a finite, non empty alphabet and $\G\subseteq W^{<\omega}(\Sigma)$. 
\begin{itemize} 
\item[(i)] Set $\left\langle \emptyset \right\rangle = \emptyset $ and for every $\bt = (t_1,\ldots,t_m) \in W^{<\omega}(\Sigma ; \upsilon)$ set \item[]$  
\left\langle \bt \right\rangle = \{(t_1(\alpha_1),\ldots,t_m(\alpha_m)) : \alpha_1,\ldots,\alpha_m \in \Sigma \}$.
\item[(ii)] $\F_{\G} = \{\bt = (t_1,\ldots,t_m)\in W^{<\omega}(\Sigma ; \upsilon) : \left\langle \bt \right\rangle \subseteq \G\}$.
\item[(ii)] $\G_* = \{\bs=(s_1,\ldots,s_k)\in W^{<\omega}(\Sigma) : \bs\in \left\langle \bt \right\rangle$ for some $\bt\in(\F_{\G})_*\}$.
\item[(iii)] $\G$ is {\em hereditary\/} if $\G_* = \G$.
\end{itemize}
\end{defn}

\begin{prop}\label{prop:finitefamily}
Let $\vec{w}=(w_n)_{n\in\nat}\in W^{<\omega}(\Sigma ; \upsilon)$ be an infinite sequence of variable words over $\Sigma$.

{\rm (i)} If $\G\subseteq RW^{<\omega}(\vec{w})$ (resp. $\F\subseteq VRW^{<\omega}(\vec{w})$) is a tree, then $\G$ (resp. $\F$) is pointwise closed if and only 
if there does not exist a reduction (resp. a variable reduction) $\vec{u}=(u_n)_{n\in\nat}$ of $\vec{w}$ such that 
$(u_1,\ldots,u_n)\in \G$ (resp. $(u_1,\ldots,u_n)\in \F)$ for all $n\in\nat$.

{\rm (ii)} If $\G\subseteq RW^{<\omega}(\vec{w})$ (resp. $\F\subseteq VRW^{<\omega}(\vec{w})$) is hereditary, then 
$\G$  (resp. $\F$) is pointwise closed if 
and only if there does not exist a variable reduction $\vec{u}$ of $\vec{w}$ such that $RW^{<\omega}(\vec{u}) \subseteq \G$ (resp. $VRW^{<\omega}(\vec{u}) \subseteq \F$). Hence, if $\G$ (resp. $\F$) is hereditary and pointwise closed, then every hereditary subfamily of $\G$ (resp. of $\F$) is also pointwise closed.

{\rm (iii)} The hereditary families $\big( W^{\xi}(\Sigma) \cap RW^{<\omega}(\vec{u}) \big)_*$, $\big( W^{\xi}(\Sigma ; \upsilon) \cap VRW^{<\omega}(\vec{u}) \big)_*$ and $\big( RW^{\xi}(\vec{w}) \cap RW^{<\omega}(\vec{u}) \big)_*$, 
$\big( VRW^{\xi}(\vec{w}) \cap VRW^{<\omega}(\vec{u}) \big)_*$ are pointwise closed for every countable ordinal $\xi\geq 0$ and $\vec{u} \in VRW^{\omega}(\vec{w})$. 
\end{prop}
\begin{proof}
(i) It follows from the relating definitions and the fact that the set $RW^{<\omega}((w_1,\ldots,w_n))$ is finite for every $n\in\nat$. 

(ii) Let $\G\subseteq R^{<\omega}(\vec{w})$ be a hereditary and not pointwise closed family. Then $\G_* = \G$, thus  $(\F_{\G})_*= \F_{\G}$. Since $\G$ is a tree, according to (i), there exists a reduction $\vec{u}=(u_n)_{n\in\nat}$ of $\vec{w}$ such that $(u_1,\ldots,u_n)\in \G=\G_*$ for all $n\in\nat$. Hence, for every $n\in\nat$ there exist $(s^n_1,\ldots,s^n_n)\in \F_{\G} \cap VRW^{<\omega}(\vec{w})$ and $\alpha^n_1,\ldots,\alpha^n_n\in \Sigma$ such that $u_i=s^n_i(\alpha^n_i)$ for every $i\leq n$. Since $\Sigma$ is finite, by a compactness argument we can find a variable reduction $\vec{s}=(s_n)_{n\in\nat}$ of $\vec{w}$ and $(\alpha_n)_{n\in\nat}\in \Sigma$ such that $(s_1,\ldots,s_n)\in \F_{\G}\cap VRW^{<\omega}(\vec{w})$ and $u_n=s_n(\alpha_n)$ for all $n\in\nat$. So, $\vec{s}\in VRW^{\omega}(\vec{w})$ and $RW^{<\omega}(\vec{s}) \subseteq \G$.

(iii) It follows from (ii).
\end{proof}

For hereditary and pointwise closed families $\G\subseteq RW^{<\omega}(\vec{w})$, $\F\subseteq VRW^{<\omega}(\vec{w})$ for some $\vec{w}\in W^{<\omega}(\Sigma ; \upsilon)$ can be defined the strong Cantor-Bendixson index $sO_{\vec{u}}(\G)$ of $\G$ and $sO_{\vec{u}}(\F)$ of $\F$ with respect to every $\vec{u} \in VRW^{\omega}(\vec{w})$.

\begin{defn}\label{def:Cantor-Bendix}
Let $\vec{w}=(w_n)_{n\in\nat} \in W^\omega (\Sigma ; \upsilon)$ be an infinite sequence of variable words over a finite, non-empty alphabet $\Sigma$ and $\G\subseteq RW^{<\omega}(\vec{w})$, $\F\subseteq VRW^{<\omega}(\vec{w})$ be hereditary and pointwise closed families. For a variable reduction $\vec{u} = (u_n)_{n\in\nat} \prec \vec{w}$ of $\vec{w}$ over $\Sigma$ we define the {\em strong Cantor-Bendixson derivatives} $(\G)_{\vec{u}}^\xi$ of $\G$, $(\F)_{\vec{u}}^\xi$ of $\F$ on $\vec{u}$ for every $\xi <\omega_1$ as follows:  

For every $\bs = (s_1,\ldots,s_k)\in \G\cap RW^{<\omega}(\vec{u})$ and $\bt = (t_1,\ldots,t_k)\in \F\cap VRW^{<\omega}(\vec{u})$ set
\begin{itemize}
\item[] $A^{\G}_{\bs}=\{w \in RW(\vec{u}) : s_1\ast\ldots\ast s_{k}\propto w, (s_1,\ldots,s_k, w-s_1\ast\ldots\ast s_{k})\notin \G \}$;    
\item[] $A^{\F}_{\bt}=\{w \in VRW(\vec{u}) : t_1\ast\ldots\ast t_{k}\propto w, (t_1,\ldots,t_k, w-t_1\ast\ldots\ast t_{k})\notin \F\}$; and
\item[] $A^{\G}_{\emptyset}=\{w\in RW(\vec{u}) :(w)\notin \G \} $, $A^{\F}_{\emptyset}=\{w\in VRW(\vec{u}) :(w)\notin \F\} $.
\end{itemize}
Then
\begin{center}
$(\G)_{\vec{u}}^0 = \{\bs \in \G\cap R^{<\omega}(\vec{u}) : A^{\G}_{\bs}$ does not contain any sequence $(w_n)_{n\in\nat}$ 
\end{center}
\begin{center}
with $w_n\propto w_{n+1}$ for every $n\in\nat$$\}$,
\end{center}
\begin{center}
$(\F)_{\vec{u}}^0 = \{\bt \in \F\cap VR^{<\omega}(\vec{u}) : A^{\F}_{\bt}$ does not contain any sequence $(w_n)_{n\in\nat}$ 
\end{center}
\begin{center}
with $w_n\propto w_{n+1}$ for every $n\in\nat$$\}$.
\end{center}
It is easy to verify that $(\G)_{\vec{u}}^0$, $(\F)_{\vec{u}}^0$ are hereditary, hence pointwise closed (Proposition~\ref{prop:finitefamily}, (ii)).
So, we can define for every $\xi >0$ the $\xi$-derivatives of $\G$ and $\F$ 
recursively as follows:

\begin{center}
 $ (\G)_{\vec{u}}^{\zeta +1} = ((\G)_{\vec{u}}^\zeta  )_{\vec{u}}^0$, $(\F)_{\vec{u}}^{\zeta +1} = ((\F)_{\vec{u}}^\zeta  )_{\vec{u}}^0$ for all $ \zeta <\omega_1$ and;
\end{center}
\begin{center}
 $(\G)_{\vec{u}}^\xi = \bigcap_{\beta<\xi} (\G)_{\vec{u}}^\beta$, $(\F)_{\vec{u}}^\xi = \bigcap_{\beta<\xi} (\F)_{\vec{u}}^\beta$ for $\xi$ a limit ordinal.
\end{center}

The {\em strong Cantor-Bendixson index} $sO_{\vec{u}}(\G)$ of $\G$ on $\vec{u}$ 
is the smallest countable ordinal $\xi$ such that $(\G)_{\vec{u}}^\xi = \emptyset$ and respectively the {\em strong Cantor-Bendixson index} $sO_{\vec{u}}(\F)$ of $\F$ on $\vec{u}$ 
is the smallest countable ordinal $\xi$ such that $(\F)_{\vec{u}}^\xi = \emptyset$.
\end{defn}

\begin{remark}\label{rem:Cantor-Bendix}
Let $\vec{w}=(w_n)_{n\in\nat} \in W^\omega (\Sigma ; \upsilon)$ and $\G\subseteq RW^{<\omega}(\vec{w})$, $\F\subseteq VRW^{<\omega}(\vec{w})$ hereditary and pointwise closed families.

{\rm (i)} The strong Cantor-Bendixson index $sO_{\vec{u}}(\G)$ and also the index $sO_{\vec{u}}(\F)$ on a variable reduction $\vec{u} \prec \vec{w}$ of $\vec{w}$ over $\Sigma$ is a countable successor ordinal less than or 
equal to the ``usual'' Cantor-Bendixson index $O(\G)$ of $\G$ and $O(\F)$ of $\F$ respectively into $\{0,1\}^{[D]^{<\omega}}$ (see \cite{K}).

{\rm (ii)} $sO_{\vec{u}} \big(\G\cap RW^{<\omega}(\vec{u})\big) = sO_{\vec{u}}(\G)$ and $sO_{\vec{u}} \big(\F\cap VRW^{<\omega}(\vec{u})\big) = sO_{\vec{u}}(\F)$.

{\rm (iii)} $sO_{\vec{u}}(\G_1) \le sO_{\vec{u}}(\G_2)$ if $\G_1,\G_2 \subseteq RW^{<\omega}(\vec{w})$ are hereditary and pointwise closed families with $\G_1\subseteq \G_2$ and also $sO_{\vec{u}}(\F_1) \le sO_{\vec{u}}(\F_2)$ if $\F_1,\F_2 \subseteq VRW^{<\omega}(\vec{w})$ are hereditary and pointwise closed families with $\F_1\subseteq \F_2$.

{\rm (iv)} If $\bs=(s_1,\ldots,s_k)\in (\G)_{\vec{u}}^\xi$ and $\vec{u}_1 \prec \vec{u}\prec \vec{w}$, then $\emptyset\in (\G)_{\vec{u}_1}^\xi$ and  $\bs_1 \in (\G)_{\vec{u}_1}^\xi$ where $\bs_1= (t_1,t_2-t_1,\ldots,t_l-t_{l-1})$ in case $\{s_1,s_1\ast s_2,\ldots,s_1\ast\ldots\ast s_k \} \cap RW(\vec{u}_1)= \{t_1,\ldots,t_l \}$, since $RW(\vec{u}_1)\subseteq RW(\vec{u})$.

{\rm (v)} If $\vec{u}_1\prec \vec{u}\prec \vec{w}$, then $sO_{\vec{u}_1} (\G) \geq sO_{\vec{u}}(\G)$ and $sO_{\vec{u}_1} (\F) \geq sO_{\vec{u}}(\F)$, according to (iv).

{\rm (vi)} Let $\vec{u}\prec \vec{w}$, $\sigma(\vec{u})=\{u_1,u_1\ast u_2,u_1\ast u_2\ast u_3,\ldots \}$ and $\vec{u}_1 \prec \vec{w}$. If  $\sigma(\vec{u}_1)\setminus \sigma(\vec{u})$ is a finite set, then $sO_{\vec{u}_1}(\G)\geq sO_{\vec{u}}(\G)$ and $sO_{\vec{u}_1}(\F)\geq sO_{\vec{u}}(\F)$.
\end{remark}
\begin{prop} 
\label{prop:Cantor-Bendix}
Let $\vec{w}=(w_n)_{n\in\nat} \in W^\omega (\Sigma ; \upsilon)$ be an infinite sequence of variable words over $\Sigma$, $\vec{u}_1=(u^1_n)_{n\in\nat}\prec \vec{u}=(u_n)_{n\in\nat}\prec \vec{w}$ be variable reductions of $\vec{w}$ over $\Sigma$ and $\xi\geq 0$ a countable ordinal. Then 

$sO_{\vec{u}_1} \big(\big(W^{\xi}(\Sigma) \cap RW^{<\omega}(\vec{u}) \big)_* \big) = sO_{\vec{u}_1} \big(\big( W^{\xi}(\Sigma ; \upsilon) \cap VRW^{<\omega}(\vec{u}) \big)_* \big) = \xi+1$, and

$sO_{\vec{u}_1} \big(\big( RW^{\xi}(\vec{w}) \cap RW^{<\omega}(\vec{u}) \big)_* \big) = sO_{\vec{u}_1} \big(\big( VRW^{\xi}(\vec{w}) \cap VRW^{<\omega}(\vec{u}) \big)_* \big) = \xi+1$.
\end{prop}

\begin{proof} 
We will prove only that $sO_{\vec{u}_1} \big(\big( VRW^{\xi}(\vec{w}) \cap VRW^{<\omega}(\vec{u}) \big)_* \big) = \xi+1\ \text{ for every }\ \xi<\omega_1\ $ and we will leave the proof of the other equalities to the reader. We mention that  $W^{\xi}(\Sigma) = RW^{\xi}(\vec{e})$ and 
$W^{\xi}(\Sigma ; \upsilon) = VRW^{\xi}(\vec{e})$ for every countable ordinal $\xi\ge 0$, in case  $\vec{e}=(e_n)_{n\in\nat}$ with $e_n=\upsilon$ for every $n\in \nat$.

For every $0<\xi<\omega_1$, the families $\big( VRW^{\xi}(\vec{w}) \cap VRW^{<\omega}(\vec{u})\big)_*$ 
are pointwise closed (Proposition~\ref{prop:finitefamily}, (iii)) and 
\begin{center}
$\big(VRW^\xi(\vec{w}) \cap VRW^{<\omega}(\vec{u})\big)(t) = VRW^{\xi_n}(\vec{w})\cap \big(VRW^{<\omega}(\vec{u}) - t \big)$  
for some $\xi_n <\xi\ ,$
\end{center} 
for every $t\in VRW(\vec{u})$ with $t\in VRW(w_1,\ldots,w_{n-1})$ for $n\in\nat$, $n>1$ (Proposition~\ref{justification 2}). 

We will prove by induction that 
$\big(\big( VRW^{\xi}(\vec{w}) \cap VRW^{<\omega}(\vec{u}) \big)_* \big)_{\vec{u}_1}^\xi = \{\emptyset\}$ for every  $0\leq\xi<\omega_1$.
Of course, 
$\big( VRW^{0}(\vec{w}) \cap VRW^{<\omega}(\vec{u}) \big)_* = \{(s) :s\in VRW(\vec{u})\} \cup \{\emptyset\}$.
Thus we have that $\big(\big( VRW^{0}(\vec{w}) \cap VRW^{<\omega}(\vec{u}) \big)_* \big)_{\vec{u}_1}^0 = \{\emptyset\}$.

Let $\xi>0$ and assume that 
$\big(\big( VRW^{\zeta}(\vec{w}) \cap VRW^{<\omega}(\vec{u}) \big)_*\big)_{\vec{u}_1}^{\zeta} = \{\emptyset\}$ for every  
$\zeta <\xi$ and $\vec{u}_1 \prec \vec{u}.$
Hence, for $\vec{u}_1 \prec \vec{u}$ and $t\in VRW(\vec{u}_1)$ with $t\in VRW(w_1,\ldots,w_{n-1})$ we have that 
$\big(\big((VRW^\xi(\vec{w}) \cap VRW^{<\omega}(\vec{u}))(t)\big)_* \big)_{\vec{u}_1}^{\xi_n} = 
\big(\big(VRW^{\xi_n}(\vec{w})\cap (VRW^{<\omega}(\vec{u}) - t)\big)_* \big)_{\vec{u}_1}^{\xi_n}= \{\emptyset\}$. 
This gives that 
$(t) \in \big(\big(VRW^\xi(\vec{w}) \cap VRW^{<\omega}(\vec{u})\big)_*\big)_{\vec{u}_1}^{\xi_n}$.
So, in case $\xi = \zeta+1$ is a successor ordinal, we have that
$(t) \in \big(\big(VRW^\xi(\vec{w}) \cap VRW^{<\omega}(\vec{u})\big)_*\big)_{\vec{u}_1}^{\zeta}$ 
for every $t\in VRW(\vec{u}_1)$, 
hence $\emptyset \in \big(\big(VRW^\xi(\vec{w}) \cap VRW^{<\omega}(\vec{u})\big)_* \big)_{\vec{u}_1}^{\xi}$. 
In case $\xi$ is a limit ordinal, we have that $\emptyset \in \big(\big(VRW^\xi(\vec{w}) \cap VRW^{<\omega}(\vec{u})\big)_* \big)_{\vec{u}_1}^{\xi}$, since 
$\emptyset\in \big(\big(VRW^\xi(\vec{w}) \cap VRW^{<\omega}(\vec{u})\big)_*\big)_{\vec{u}_1}^{\xi_n}$ for every $n\in\nat$ 
and $\sup \xi_n  =\xi$.

If $\{\emptyset\} \ne \big(\big(VRW^\xi(\vec{w}) \cap VRW^{<\omega}(\vec{u})\big)_* \big)_{\vec{u}_1}^{\xi}$ 
for some $\vec{u}_1 \prec \vec{u}$, 
then there exist $\vec{u}_2 \prec \vec{u}_1$ and  $t\in VRW(\vec{u}_2)$ such that
$\big(\big((VRW^\xi(\vec{w}) \cap VRW^{<\omega}(\vec{u}))(t)\big)_* \big)_{\vec{u}_2}^{\xi} = 
\big(\big(VRW^{\xi_n}(\vec{w})\cap (VRW^{<\omega}(\vec{u}) - t)\big)_*\big)_{\vec{u}_2}^{\xi}\ne \{\emptyset\}$ 
(see Lemma 2.8 in \cite{F3}). 
This is a contradiction to the induction hypothesis. 
Hence, $\{\emptyset\} \ne \big(\big(VRW^\xi(\vec{w}) \cap VRW^{<\omega}(\vec{u})\big)_* \big)_{\vec{u}_1}^{\xi}$ 
and $sO_{\vec{u}_1} \big(\big( VRW^{\xi}(\vec{w}) \cap VRW^{<\omega}(\vec{u}) \big)_* \big) = \xi+1$ 
for every $\xi<\omega_1$ and $\vec{u}_1 \prec \vec{u}\prec \vec{w}$ .
\end{proof}

\begin{cor}\label{cor:tree} 
For every $\vec{w}=(w_n)_{n\in\nat} \in W^\omega (\Sigma ; \upsilon)$ and countable ordinals $\xi_1,\xi_2$ with $\xi_1<\xi_2$ there exists a variable reduction $\vec{u}\prec \vec{w}$ of $\vec{w}$ over $\Sigma$ such that:
\begin{center}
$\big(W^{\xi_1}(\Sigma)\big)_* \cap RW^{<\omega}(\vec{u}) \subseteq \big(W^{\xi_2}(\Sigma)\big)^* \setminus W^{\xi_2}(\Sigma)$ and
\end{center}
\begin{center}
$\big(W^{\xi_1}(\Sigma ; \upsilon)\big)_* \cap VRW^{<\omega}(\vec{u}) \subseteq \big(W^{\xi_2}(\Sigma ; \upsilon)\big)^* \setminus  W^{\xi_2}(\Sigma ; \upsilon)$.
\end{center}
\end{cor}

\begin{proof} 
Of course $\big(W^{\xi_1}(\Sigma)\big)_*\subseteq W^{<\omega}(\Sigma)$ and $\big(W^{\xi_1}(\Sigma ; \upsilon)\big)_*\subseteq  W^{<\omega}(\Sigma ; \upsilon)$ are trees. According to Theorem~\ref{thm:treeRam}, for every infinite sequence $\vec{w} \in W^\omega (\Sigma ; \upsilon)$ there exists a variable reduction $\vec{u} \prec \vec{w} $ of $\vec{w}$ over $\Sigma$ such that :

either $W^{\xi_2}(\Sigma)\cap RW^{<\omega}(\vec{u})\subseteq \big(W^{\xi_1}(\Sigma)\big)_*$ 

or $\big(W^{\xi_1}(\Sigma)\big)_* \cap RW^{<\omega}(\vec{u}) \subseteq \big(W^{\xi_2}(\Sigma)\big)^* \setminus W^{\xi_2}(\Sigma)$, and

either $ W^{\xi_2}(\Sigma ; \upsilon)\cap VRW^{<\omega}(\vec{u})\subseteq \big(W^{\xi_1}(\Sigma ; \upsilon)\big)_*$ 

or $\big(W^{\xi_1}(\Sigma ; \upsilon)\big)_*\cap VRW^{<\omega}(\vec{u}) \subseteq \big(W^{\xi_2}(\Sigma ; \upsilon)\big)^* \setminus W^{\xi_2}(\Sigma ; \upsilon)$.

The first alternative in each of the two dichotomies is impossible, since, otherwise 
according to Proposition~\ref{prop:Cantor-Bendix}, $\xi_2 +1 = sO_{\vec{u}} \big(\big(W^{\xi_2}(\Sigma) \cap RW^{<\omega}(\vec{u}) \big)_* \big) \le 
sO_{\vec{u}} \big( \big(W^{\xi_1}(\Sigma)\big)_* \big) = \xi_1 +1$ or 
$\xi_2 +1 = sO_{\vec{u}} \big(\big(W^{\xi_2}(\Sigma ; \upsilon) \cap VRW^{<\omega}(\vec{u}) \big)_* \Big) \le 
sO_{\vec{u}} \big(\big(W^{\xi_1}(\Sigma ; \upsilon)\big)_* \big) = \xi_1 +1$; a contradiction. 
\end{proof}

\section{Schreier-type extension of Carlson's Nash-Williams type partition theorem for words}

According to the partition theorem on Schreier families proved in Section 2, for every countable ordinal $\xi$, every non-empty, finite alphabet $\Sigma$  and 
every partition $\G$ of the set $W^{<\omega}(\Sigma)$ of all the finite sequences of words over $\Sigma$, 
there exists an infinite sequence $\vec{u}$ of variable words over $\Sigma$ , 
all of whose finite reductions in the Schreier family $W^{\xi}(\Sigma)$ are either in the partition family $\G$ itself or in the complement $W^{<\omega}(\Sigma)\setminus \G$; but Theorem~\ref{xi Carlson} 
naturally can provide no information whatsoever on whether all these finite reductions are in $\G$ or in its complement $W^{<\omega}(\Sigma)\setminus \G$. 
In this section we will obtain, 
for a partition family $\G$ that is a tree, a criterion on this matter, in 
terms of the strong Cantor-Bendixson index of $\G$: if 
this index is greater than $\xi +1$, all $W^{\xi}(\Sigma)$-finite reductions fall in $\G$, 
and if less than $\xi$, in $W^{<\omega}(\Sigma)\setminus \G$ (albeit in a weaker, non-symmetrical 
manner) (Theorem~\ref{block-NashWilliams1} and Theorem~\ref{block-NashWilliams2}).

It will be observed that the main dichotomy of Theorem~\ref{block-NashWilliams1} 
is non-symmetric, reflecting the fact that the treeness property is assumed for the family $\G$ itself only, and of course not for its complement $W^{<\omega}(\Sigma)\setminus \G$. This type of 
non-symmetric dichotomies is characteristic of Nash-Williams type partition theorem; in fact, from Theorem~\ref{block-NashWilliams1} and the analogous Theorem~\ref{block-NashWilliams2} for variable words, we 
will derive in the sequel various strong forms of Nash-Williams type partition theorems for words and variable words involving the Schreier-type families of words and the Cantor-Bendixson index (Theorem B, Corollaries~\ref{cor:gowers}, \ref{cor:gowers3}, \ref{cor:gowers1},  \ref{cor:blockNW1} ), which imply as well as Carlson's infinitary partition theorem (Theorem~\ref{cor:IblockNW}, \cite{C}). 

In the proof of Theorem~\ref{block-NashWilliams1} bolow we use Theorem~\ref{xi Carlson} and also 
we exploit the properties of the Schreier-type families $\W^{\xi}(\Sigma)$ for $\xi<\omega_1$ proved in the previous Section 3. Towards this purpose we introduce the following definition. 

\begin{defn}\label{def:blockNW} 
Let $\G\subseteq W^{<\omega}(\Sigma)$ and $\F\subseteq W^{<\omega}(\Sigma ; \upsilon)$. We set 
\begin{itemize} 
\item[(i)] $\G_0 = \{\bs \in \G : \bs\in \left\langle \bt \right\rangle$ for some $\bt\in\F_{\G}\}$.
\item[(ii)] $\G_h = \{\bs \in \G_0 $ : in case $\bs\in \left\langle \bt \right\rangle$ for some $\bt\in\F_{\G}$ then $\left\langle \bu \right\rangle\subseteq\G$ 
for every $\bu\in VRW^{<\omega}(\bt_1)$ for $\bt_1\propto \bt \} 
\cup\{\emptyset\}$.

\item[(iii)] $\F_h = \{\bt \in \F: VRW^{<\omega}(\bt_1)\subseteq\F$ for every $\bt_1\propto \bt \}\cup \{\emptyset\}$.
\end{itemize}
Of course, $\G_h$, $\F_h$ are the largest subfamilies of $\G\cup \{\emptyset\}$, $\F\cup \{\emptyset\}$ which are hereditary.
\end{defn}

\begin{thm}
\label{block-NashWilliams1}
Let $\G\subseteq W^{<\omega}(\Sigma)$ be a family of finite sequences of words over the finite, non-empty alphabet $\Sigma$ which is a tree and $\vec{w} \in W^\omega (\Sigma ; \upsilon)$ be an infinite sequence of variable words over $\Sigma$. We have the following cases:

\noindent {\bf [Case 1]}
The family $\G_h \cap RW^{<\omega}(\vec{w})$ is not pointwise closed. 

Then, there exists a variable reduction $\vec{u}$ of $\vec{w}$ over $\Sigma$ such that 
\begin{center}
$RW^{<\omega}(\vec{u})\subseteq \G$.
\end{center}
\noindent {\bf [Case 2]}
The family $\G_h \cap RW^{<\omega}(\vec{w})$ is pointwise closed. 

Then, setting
\begin{equation*}
\zeta_{\vec{w}}^{\G} = \xi_{\vec{w}}^{\G_h} = \sup \{sO_{\vec{u}} \big(\G_h\cap RW^{<\omega}(\vec{w})\big) : \vec{u} \prec \vec{w}\}\ ,
\end{equation*}
which is a countable ordinal, the following subcases obtain:
\begin{itemize}
\item[2(i)] If $\xi+1 <\zeta_{\vec{w}}^{\G}$, then there exists $\vec{u} \prec \vec{w}$ 
such that 
\begin{center}
$W^{\xi}(\Sigma)\cap RW^{<\omega}(\vec{u})\subseteq \G$; 
\end{center}
\item[2(ii)] if $\omega_1 > \xi+1>\xi>\zeta_{\vec{w}}^{\G}$, then for every $\vec{u} \prec \vec{w}$ there exists 
$\vec{u}_1 \prec \vec{u}$ such that 
\begin{center}
$W^{\xi}(\Sigma)\cap RW^{<\omega}(\vec{u}) \subseteq W^{<\omega}(\Sigma)\setminus \G$
\end{center}
(equivalently $\G\cap RW^{<\omega}(\vec{u}) \subseteq \big(W^{\xi}(\Sigma)\big)^* \setminus W^{\xi}(\Sigma)$) ; and 
\item[2(iii)] if $\xi+1 = \zeta_{\vec{s}_0}^{\G}$ or $\xi = \zeta_{\vec{s}_0}^{\G}$, then there exists $\vec{u} \prec \vec{w}$ 
such that 
\item[{}] either $W^{\xi}(\Sigma)\cap RW^{<\omega}(\vec{u})\subseteq \G$ or 
$W^{\xi}(\Sigma)\cap RW^{<\omega}(\vec{u}) \subseteq W^{<\omega}(\Sigma)\setminus \G$.
\end{itemize}
\end{thm}

\begin{proof}{} 
[Case 1] If the hereditary family $\G_h \cap RW^{<\omega}(\vec{w})$ is not pointwise 
closed, then, according to Proposition~\ref{prop:finitefamily}, there exists $\vec{u} \prec \vec{w}$ such that 
\begin{itemize}
\item[{}] $RW^{<\omega}(\vec{u})\subseteq \G_h\cap RW^{<\omega}(\vec{w})\subseteq \G_h \subseteq \G$.
\end{itemize}

\noindent [Case 2] 
If the hereditary family $\G_h \cap RW^{<\omega}(\vec{w})$ is pointwise closed, 
then the index $\zeta_{\vec{w}}^{\G}=\xi_{\vec{w}}^{\G_h}$ is countable, since the ``usual'' Cantor-Bendixson 
index $O(\G_h \cap RW^{<\omega}(\vec{w}))$ of $\G_h\cap RW^{<\omega}(\vec{w})$ into $\{0,1\}^{[D]^{<\omega}}$ is countable and $sO_{\vec{u}}\big(\G_h\cap RW^{<\omega}(\vec{w})\big)\le O \big(\G_h\cap RW^{<\omega}(\vec{w})\big)$ for every $\vec{u} \prec \vec{w}$ (Remark~\ref{rem:Cantor-Bendix}(i)and (ii)).

2(i) Let $\xi+1 < \zeta_{\vec{w}}^{\G}$. Then $\xi+1 < \xi_{\vec{w}}^{\G_h}$, so there exists $\vec{u}_1\prec\vec{w}$ such that $\xi+1 < sO_{\vec{u}_1}\big(\G_h \cap RW^{<\omega}(\vec{w})\big)$. 
According to Theorem~\ref{thm:treeRam}, there exists a variable reduction $\vec{u}\prec \vec{u}_1$ of $\vec{u}_1$ over $\Sigma$ such that  

either $W^{\xi}(\Sigma)\cap RW^{<\omega}(\vec{u})\subseteq \G_h$ or 
$\G_h\cap RW^{<\omega}(\vec{u}) \subseteq \big(W^{\xi}(\Sigma)\big)^* \setminus W^{\xi}(\Sigma)\subseteq \big(W^{\xi}(\Sigma)\big)^*$. 
\newline
The second alternative is impossible. 
Indeed, if $\G_h\cap RW^{<\omega}(\vec{u}) \subseteq \big(W^{\xi}(\Sigma)\big)^*$, then, according 
to Remark ~\ref{rem:Cantor-Bendix} and Proposition~\ref{prop:Cantor-Bendix}, 
\newline
$sO_{\vec{u}_1}\big(\G_h \cap RW^{<\omega}(\vec{w})\big) \leq sO_{\vec{u}} \big(\G_h\cap RW^{<\omega}(\vec{w})\big) = sO_{\vec{u}}\big(\G_h \cap RW^{<\omega}(\vec{u})\big) \leq sO_{\vec{u}}\big(\big(W^{\xi}(\Sigma)\big)^* \big) = \xi+1$; 
a contradiction. Hence, $W^{\xi}(\Sigma)\cap RW^{<\omega}(\vec{u})\subseteq \G_h\subseteq \G$.

2(ii) Let $\xi +1>\xi> \zeta_{\vec{w}}^{\G}$, and $\vec{u}\prec\vec{w}$. For every countable ordinal $\zeta$ with $\zeta +1 > \zeta_{\vec{w}}^{\G}$ there exists 
a variable reduction $\vec{u}_1\prec \vec{u}$ of $\vec{u}$ over $\Sigma$ such that  
\begin{itemize}
\item[{}] $W^{\zeta}(\Sigma)\cap RW^{<\omega}(\vec{u}_1) \subseteq W^{<\omega}(\Sigma)\setminus \G_h$,
\end{itemize}
Indeed, according to the partition theorem on Schreier families (Theorem~\ref{xi Carlson}), there exists 
a variable reduction $\vec{u}_1\prec \vec{u}$ of $\vec{u}$ over $\Sigma$ such that  
\begin{itemize}
\item[{}] either $W^{\zeta}(\Sigma)\cap RW^{<\omega}(\vec{u}_1)\subseteq \G_h$ or 
$W^{\zeta}(\Sigma)\cap RW^{<\omega}(\vec{u}_1) \subseteq W^{<\omega}(\Sigma)\setminus \G_h$,
\end{itemize}
The first alternative is impossible, since if $W^{\zeta}(\Sigma)\cap RW^{<\omega}(\vec{u}_1)\subseteq \G_h$, then, according to Remark ~\ref{rem:Cantor-Bendix} and Proposition~\ref{prop:Cantor-Bendix}, we obtain that 
\begin{equation*}
\zeta+1 = sO_{\vec{u}_1} \big(\big(W^{\zeta}(\Sigma)\cap RW^{<\omega}(\vec{u}_1)\big)_* \big) \le sO_{\vec{u}_1} 
\big(\G_h \cap RW^{<\omega}(\vec{u}_1)\big) \le \xi_{\vec{w}}^{\G_h}=\zeta_{\vec{w}}^{\G}\ ;
\end{equation*}
a contradiction. Hence, there exists a variable reduction $\vec{u}_1\prec \vec{u}$ of $\vec{u}$ over $\Sigma$ such that
\begin{itemize}
\item[{}] $W^{\zeta_{\vec{w}}^{\G}}(\Sigma)\cap RW^{<\omega}(\vec{u}_1) \subseteq W^{<\omega}(\Sigma)\setminus \G_h$.
\end{itemize}

According to Theorem~\ref{xi Carlson},
there exists $\vec{u}_2 \prec \vec{u}_1$ such that 
\begin{itemize}
\item[{}] either $W^{\xi}(\Sigma)\cap RW^{<\omega}(\vec{u}_2)\subseteq \G$ or 
$W^{\xi}(\Sigma)\cap RW^{<\omega}(\vec{u}_2) \subseteq W^{<\omega}(\Sigma)\setminus \G$.
\end{itemize}
We claim that the first alternative does not hold. 
Indeed, if $W^{\xi}(\Sigma)\cap RW^{<\omega}(\vec{u}_2)\subseteq \G$, then 
$\big(W^{\xi}(\Sigma)\cap RW^{<\omega}(\vec{u}_2)\big)^* \subseteq \G^* = \G$. 
Using the canonical representation of every infinite sequence of words over $\Sigma$ with respect to $W^{\xi}(\Sigma)$ (Proposition~\ref{prop:canonicalrep}) we have that 
\begin{itemize}
\item[{}] $\big(W^{\xi}(\Sigma)\big)^* \cap RW^{<\omega}(\vec{u}_2) = \big(W^{\xi}(\Sigma)\cap RW^{<\omega}(\vec{u}_2)\big)^*$.
\end{itemize}
Hence, $\big(W^{\xi}(\Sigma)\big)^* \cap RW^{<\omega}(\vec{u}_2) \subseteq \G$.
\newline
Since $\xi >\zeta_{\vec{w}}^{\G}$, according to Corollary~\ref{cor:tree}, there exists $\vec{u}_3 \prec\vec{u}_2$ such that 
\begin{itemize}
\item[{}] $\big(W^{\zeta_{\vec{w}}^{\G}}(\Sigma)\big)_* \cap RW^{<\omega}(\vec{u}_3)\subseteq \big(W^{\xi}(\Sigma)\big)^* \cap RW^{<\omega}(\vec{u}_2) \subseteq \G$.
\end{itemize}
Thus $\big(W^{\zeta_{\vec{w}}^{\G}}(\Sigma)\big)_* \cap RW^{<\omega}(\vec{u}_3) \subseteq \G_h$, since $\big(W^{\zeta_{\vec{w}}^{\G}}(\Sigma)\big)_* \cap RW^{<\omega}(\vec{u}_3)$ is a hereditary family. 
This is a contradiction, since $\vec{u}_3\prec\vec{u}_1$ and $W^{\zeta_{\vec{w}}^{\G}}(\Sigma)\cap RW^{<\omega}(\vec{u}_1) \subseteq W^{<\omega}(\Sigma)\setminus \G_h$. Hence, 
\begin{itemize}
\item[{}] $W^{\xi}(\Sigma)\cap RW^{<\omega}(\vec{u}_2) \subseteq W^{<\omega}(\Sigma)\setminus \G$ and 
$\G \cap RW^{<\omega}(\vec{u}_2) \subseteq \big(W^{\xi}(\Sigma)\big)^* \setminus W^{\xi}(\Sigma)$, 
\end{itemize}
according to Proposition~\ref{prop:tree}. 

2(iii) In the cases $\xi+1 = \zeta_{\vec{s}_0}^{\G}$ or $\xi = \zeta_{\vec{s}_0}^{\G}$, use Theorem~\ref{xi Carlson}.
\end{proof}

\begin{remark}
\label{since}
Let $\G\subseteq W^{<\omega}(\Sigma)$ be a tree and let $\vec{w} \in W^\omega (\Sigma ; \upsilon)$.

{\rm (i)} That both alternatives may materialize in case $\xi +1=\zeta_{\vec{w}}^{\G}$ can be seen by considering two 
simple examples: 

(1) Let $\F = \{\bt = (t_1< t_2< \cdots < t_{2k+2}) \in W^{<\omega}(\Sigma ; \upsilon) : 
k\in\nat$ and $\min d(\bt) =k\}$ and $\G = \{\bs\in W^{<\omega}(\Sigma) :  \bs\in \left\langle \bt \right\rangle$ for some $\bt\in\F_*\}$. 
It is easy to see that the hereditary family $\G$ is pointwise closed (according to 
Proposition~\ref{prop:finitefamily}). Analogously to Proposition~\ref{prop:Cantor-Bendix}, can be proved that
$sO_{\vec{u}}\big(G\cap RW^{<\omega}(\vec{w})\big) = \omega+1$ for every $\vec{w}\in W^\omega (\Sigma ; \upsilon)$ and $\vec{u} \prec \vec{w}$.
Thus, $\zeta_{\vec{w}}^{\G}=\xi_{\vec{w}}^{\G} = \omega+1$. It is now easy to verify that 
\begin{center}
$W^{\omega}(\Sigma)\cap RW^{<\omega}(\vec{u})\subseteq \G$ for every $\vec{u} \prec \vec{w}$.
\end{center}

(2) Let $\F = \{\bt = (t_1< t_2< \cdots < t_{k+1}) \in W^{<\omega}(\Sigma ; \upsilon) : k\in\nat$ and $\min d(\bt) =2k\}$ and $\G = \{\bs\in W^{<\omega}(\Sigma) : \bs\in \left\langle \bt \right\rangle$ for some $\bt\in\F_*\}$. The hereditary family $\G$ is pointwise closed. Setting $\vec{w}=(w_n)_{n\in\nat} \in W^\omega (\Sigma ; \upsilon)$ with $w_{1}=\upsilon$ and $w_{n}=\upsilon\ast\upsilon$ for every $1<n\in\nat$, we have that $sO_{\vec{u}}\big(\G\cap RW^{<\omega}(\vec{w})\big) = \omega+1$ for every $\vec{u} \prec \vec{w}$. Thus, $\zeta_{\vec{w}}^{\G}=\xi_{\vec{w}}^{\G} = \omega+1$. 
It is now easy to see that 
\begin{center}
$W^{\omega}(\Sigma)\cap RW^{<\omega}(\vec{u}) \subseteq W^{<\omega}(\Sigma)\setminus \G$ for every $\vec{u} \prec  \vec{w}$, 
\end{center}
since
$\G\cap RW^{<\omega}(\vec{u}) \subseteq \big(W^{\omega}(\Sigma)\big)^* \setminus W^{\omega}(\Sigma)$. 

{\rm (ii)}  In case the family $\G\subseteq W^{<\omega}(\Sigma)$ is hereditary and $\xi =\zeta_{\vec{w}}^{\G} = \xi_{\vec{w}}^{\G}<\omega_1$, then can be proved that for every $\vec{u} \prec \vec{w}$ there exists 
$\vec{u}_1 \prec \vec{u}$ such that 
\begin{center}
$W^{\xi}(\Sigma)\cap RW^{<\omega}(\vec{u}) \subseteq W^{<\omega}(\Sigma)\setminus \G$.
\end{center}
\end{remark}

For a partition of all the finite sequences of variable words over $\Sigma$ which is a tree holds an analogous  strengthened theorem, which in fact is a stronger form of Carlson's infinitary partition theorem (Theorem~\ref{cor:IblockNW}, \cite{C}). Although, the proof of this theorem is analogous to the proof of Theorem~\ref{block-NashWilliams1}, for completeness we will give a sketch of it.

\begin{thm}
\label{block-NashWilliams2}
Let $\F\subseteq W^{<\omega}(\Sigma ; \upsilon)$ be a family of finite sequences of variable words over the finite, non-empty alphabet $\Sigma$ which is a tree and $\vec{w} \in W^\omega (\Sigma ; \upsilon)$ be an infinite sequence of variable words over $\Sigma$. We have the following cases: 

\noindent {\bf [Case 1]}
The family $\F_h \cap VRW^{<\omega}(\vec{w})$ is not pointwise closed. 

Then, there exists a variable reduction $\vec{u}$ of $\vec{w}$ over $\Sigma$ such that 
\begin{center}
$VRW^{<\omega}(\vec{u})\subseteq \F$.
\end{center}
\noindent {\bf [Case 2]}
The family $\F_h \cap VRW^{<\omega}(\vec{w})$ is pointwise closed. 

Then, setting
\begin{equation*}
\zeta_{\vec{w}}^{\F} = \xi_{\vec{w}}^{\F_h} = \sup \{sO_{\vec{u}} \big(\F_h\cap VRW^{<\omega}(\vec{w})\big) : \vec{u}\prec\vec{w}\}\ ,
\end{equation*}
which is a countable ordinal, the following subcases obtain:
\begin{itemize}
\item[2(i)] If $\xi+1 <\zeta_{\vec{w}}^{\F}$, then there exists $\vec{u}\prec\vec{w}$ 
such that 
\begin{center}
$W^{\xi}(\Sigma ; \upsilon)\cap VRW^{<\omega}(\vec{u})\subseteq \F$;
\end{center}
\item[2(ii)] if $\xi+1>\xi>\zeta_{\vec{w}}^{\F}$, then for every $\vec{u} \prec \vec{w}$ there exists 
$\vec{u}_1 \prec \vec{u}$ such that 
\begin{center}
$W^{\xi}(\Sigma ; \upsilon)\cap VRW^{<\omega}(\vec{u}) \subseteq W^{<\omega}(\Sigma ; \upsilon)\setminus \F$
\end{center}
(equivalently $\F\cap VRW^{<\omega}(\vec{u}) \subseteq \big(W^{\xi}(\Sigma ; \upsilon)\big)^* \setminus W^{\xi}(\Sigma ; \upsilon)$) ; and 
\item[2(iii)] if $\xi+1 = \zeta_{\vec{s}_0}^{\F}$ or $\xi = \zeta_{\vec{s}_0}^{\F}$, then there exists $\vec{u} \prec \vec{w}$ 
such that 
\item[{}] either $W^{\xi}(\Sigma ; \upsilon)\cap VRW^{<\omega}(\vec{u})\subseteq \F$ or 
$W^{\xi}(\Sigma ; \upsilon)\cap VRW^{<\omega}(\vec{u}) \subseteq W^{<\omega}(\Sigma ; \upsilon)\setminus \F$.
\end{itemize}
\end{thm}
\begin{proof}{} 
[Case 1] If the hereditary family $\F_h \cap VRW^{<\omega}(\vec{w})$ is not pointwise 
closed, then, there exists $\vec{u}\prec \vec{w}$ with $VRW^{<\omega}(\vec{u}) \subseteq \F_h\cap VRW^{<\omega}(\vec{w})\subseteq \F$ (Proposition~\ref{prop:finitefamily}). 

\noindent [Case 2] If the hereditary family $\F_h \cap VRW^{<\omega}(\vec{w})$ is pointwise closed, 
then the index $\xi_{\vec{w}}^{\F_h}$ is countable, according to Remark~\ref{rem:Cantor-Bendix}(i)and (ii). 

2(i) Let $\xi+1 <\zeta_{\vec{w}}^{\F}$. Then there exists $\vec{u}_1 \prec \vec{w}$ such that $\xi+1 < sO_{\vec{u}_1}\big(\F_h \cap VRW^{<\omega}(\vec{w})\big)$. Using Theorem~\ref{thm:treeRam}, Remark ~\ref{rem:Cantor-Bendix} and Proposition~\ref{prop:Cantor-Bendix}, we have that 
\begin{itemize}
\item[{}]
$W^{\xi}(\Sigma ; \upsilon)\cap VRW^{<\omega}(\vec{u})\subseteq \F_h\subseteq \F$ .
\end{itemize}
2(ii) Let $\xi+1>\xi>\zeta_{\vec{w}}^{\F}$ and $\vec{u} \prec \vec{w}$. 
According to Theorem~\ref{xi Carlson}, Remark ~\ref{rem:Cantor-Bendix} and Proposition~\ref{prop:Cantor-Bendix}, there exists a variable reduction $\vec{u}_1\prec \vec{u}$ of $\vec{u}$ over $\Sigma$ such that  
\begin{itemize}
\item[{}] $W^{\zeta_{\vec{w}}^{\F}}(\Sigma ; \upsilon)\cap VRW^{<\omega}(\vec{u}_1) \subseteq W^{<\omega}(\Sigma ; \upsilon)\setminus \F_h$.
\end{itemize}
Using again Theorem~\ref{xi Carlson}, there exists $\vec{u}_2 \prec \vec{u}_1$ such that 
\begin{itemize}
\item[{}] either $W^{\xi}(\Sigma ; \upsilon)\cap VRW^{<\omega}(\vec{u}_2)\subseteq \F$ or 
$W^{\xi}(\Sigma)\cap VRW^{<\omega}(\vec{u}_2) \subseteq W^{<\omega}(\Sigma)\setminus \F$.
\end{itemize}
We claim that the first alternative does not hold. 
Indeed, if $W^{\xi}(\Sigma ; \upsilon)\cap VRW^{<\omega}(\vec{u}_2)\subseteq \F$, then using the canonical representation of every infinite sequence of variable words over $\Sigma$ with respect to $W^{\xi}(\Sigma ; \upsilon)$ (Proposition~\ref{prop:canonicalrep}) it is easy to check that 
\begin{itemize}
\item[{}] $\big(W^{\xi}(\Sigma ; \upsilon)\big)^* \cap VRW^{<\omega}(\vec{u}_2) = \big(W^{\xi}(\Sigma ; \upsilon)\cap VRW^{<\omega}(\vec{u}_2)\big)^*\subseteq \F^* = \F$.
\end{itemize}
Since $\xi >\zeta_{\vec{w}}^{\F}$, according to Corollary~\ref{cor:tree}, there exists $\vec{u}_3 \prec\vec{u}_2$ such that 
\begin{center}
$\big(W^{\zeta_{\vec{w}}^{\F}}(\Sigma ; \upsilon)\big)_* \cap VRW^{<\omega}(\vec{u}_3)\subseteq \big(W^{\xi}(\Sigma ; \upsilon)\big)^* \cap VRW^{<\omega}(\vec{u}_2) \subseteq \F$;
\end{center}
and consequently such that $\big(W^{\zeta_{\vec{w}}^{\F}}(\Sigma ; \upsilon)\big)_* \cap VRW^{<\omega}(\vec{u}_3)\subseteq \F_h$.
This is a contradiction. 

2(iii) In the cases $\xi+1 = \zeta_{\vec{w}}^{\F}$ or $\xi = \zeta_{\vec{w}}^{\F}$, use Theorem~\ref{xi Carlson}.
\end{proof} 

That both alternatives may materialize in case $\xi +1=\zeta_{\vec{w}}^{\F}$ can be seen by considering the following 
examples:

(1) Let $\F = \{\bt = (t_1< t_2< \cdots < t_{2k+2}) \in W^{<\omega}(\Sigma ; \upsilon) : 
k\in\nat$ and $\min d(\bt) =k\}$. The hereditary family $\F_*$ is pointwise closed and 
$sO_{\vec{u}}(\F_*) = \omega+1$ for every $\vec{w}\in\W^\omega (\Sigma ; \upsilon)$ and $\vec{u} \prec \vec{w}$.
Thus, $\zeta_{\vec{w}}^{\F_*}=\xi_{\vec{w}}^{\F_*} = \omega+1$. It is now easy to verify that 
\begin{center}
$W^{\omega}(\Sigma ; \upsilon)\cap VRW^{<\omega}(\vec{u})\subseteq \F_*$ for every $\vec{u} \prec \vec{w}$.
\end{center}

(2) Let $\F = \{\bt = (t_1< t_2< \cdots < t_{k+1}) \in W^{<\omega}(\Sigma ; \upsilon) : k\in\nat$ and $\min d(\bt) =2k\}$. The hereditary family $\F_*$ is pointwise closed and $sO_{\vec{u}}(\F_*) = \omega+1$ for 
every $\vec{u} \prec \vec{w}$, where $\vec{w}=(w_n)_{n\in\nat} \in \W^\omega (\Sigma ; \upsilon)$ with $w_{1}=\upsilon$ and $w_{n}=\upsilon\ast\upsilon$ for every $1<n\in\nat$. Thus, $\zeta_{\vec{w}}^{\F_*}=\xi_{\vec{w}}^{\F_*} = \omega+1$. It is now easy to see that 
\begin{center}
$\F_* \cap VRW^{<\omega}(\vec{u}) \subseteq \big(W^{\omega}(\Sigma ; \upsilon)\big)^* \setminus W^{\omega}(\Sigma ; \upsilon)$ 
for every $\vec{u} \prec  \vec{w}$. 
\end{center}

An immediate consequence of Theorems~\ref{block-NashWilliams1} and \ref{block-NashWilliams2} is Theorem B, referred to  in the introduction, which is a strengthened form of Theorem~\ref{xi Carlson} in that the partitions are trees. A 
quite simplified consequence of Theorem~\ref{block-NashWilliams2}, one not involving Schreier-type families of countable ordinal index, is equivalent to Carlson's infinitary partition theorem (Theorem~\ref{cor:IblockNW}) proved in \cite{C}. 

\begin{cor}
\label{cor:gowers}
Let $\F \subseteq W^{<\omega}(\Sigma ; \upsilon)$ be a family of finite sequences of variable words 
over an alphabet $\Sigma$ which is a tree. Then for every infinite sequence $\vec{w} \in W^\omega (\Sigma ; \upsilon)$ of variable words over $\Sigma$ there exists a variable reduction $\vec{u}\prec\vec{w}$ of $\vec{w}$ over $\Sigma$ such that:

 either $VRW^{<\omega}(\vec{u})\subseteq \F$, 

 or for every variable reduction $\vec{u}_1 $ of $\vec{u}$ there exists an initial segment of $\vec{u}_1$ 
which belongs to $W^{<\omega}(\Sigma ; \upsilon)\setminus \F$. 

\end{cor}
\begin{proof}
The proof follows from 
Theorem~\ref{block-NashWilliams2} (case 1 and subcase 2(ii)) and Proposition~\ref{prop:canonicalrep2}. 
\end{proof}

\begin{remark}\label{rem:block-NashWilliams} 
Carlson's Theorem~\ref{cor:IblockNW} is equivalent to Corollary~\ref{cor:gowers}. Indeed: 

(i) Corollary~\ref{cor:gowers} implies Theorem~\ref{cor:IblockNW}. Indeed, let $\U \subseteq W^{\omega}(\Sigma ; \upsilon)$ be a pointwise closed family of infinite sequences of variable words over $\Sigma$ and $\vec{w} \in W^\omega (\Sigma ; \upsilon)$. Set 

$\F_{\U} = \{\bt=(t_1,\ldots,t_k) \in W^{<\omega}(\Sigma ; \upsilon)$: $k\in\nat$ and there exists $\vec{t}\in \U$ with $\bt\propto \vec{t} \}$.
\newline 
Since the family $\F_{\U}$ is a tree, we use Corollary~\ref{cor:gowers}.
Then we have the following two cases:

\noindent {[Case 1]}
There exists $\vec{u}\prec\vec{w}$ such that $VRW^{<\omega}(\vec{u})\subseteq \F_{\U}$. 
Then, $VRW^{\omega}(\vec{u})\subseteq \U$. 
Indeed, if $\vec{t} = (t_n)_{n\in\nat} \in VRW^{\omega}(\vec{u})$, then 
$(t_1,\ldots,t_k)\in \F_{\U}$ for every $k\in\nat$. 
Hence, for each $k\in\nat$ there exists $\vec{t}_k=(t^k_n)_{n\in\nat}\in\U$ such that 
$t^k_n=t_n$ for every $n\in\nat$ with $n\leq k$. 
Since $(\vec{t}_k)_{k\in\nat}$ converges 
pointwise to $\vec{t}$ and $\U$ is pointwise closed, 
we have that $\vec{t}\in\U$ and consequently that $VRW^{\omega}(\vec{u})\subseteq \U$.

\noindent {[Case 2]}
There exists $\vec{u}\prec\vec{w}$ such that every variable reduction $\vec{u}_1 $ of $\vec{u}$ has an initial segment  belonging to $W^{<\omega}(\Sigma ; \upsilon)\setminus \F_{\U}$. 
Hence, $VRW^{\omega}(\vec{u})\subseteq W^{\omega}(\Sigma ; \upsilon) \setminus \U$.

(ii) Theorem~\ref{cor:IblockNW} implies Corollary~\ref{cor:gowers}. Indeed, let $\F \subseteq W^{<\omega}(\Sigma ; \upsilon)$ which is a tree and $\vec{w} \in W^\omega (\Sigma ; \upsilon)$. Set  

$\U_{\F} = \{ \vec{t}=(t_n)_{n\in\nat}\in W^{\omega}(\Sigma ; \upsilon):$ there exists $k\in\nat$ such that  $(t_1,\ldots,t_k)\in \F \}.$
\newline
Then $\W^{\omega}(\Sigma ; \upsilon)\setminus \U_{\F}$ is pointwise closed, 
so, using Theorem~\ref{cor:IblockNW} for the family 
$\W^{\omega}(\Sigma ; \upsilon)\setminus \U_{\F}$, we obtain Corollary~\ref{cor:gowers}. 
\end{remark}

In fact Corollary~\ref{cor:gowers} holds for arbitrary partitions of $W^{<\omega}(\Sigma ; \upsilon)$, not necessarily trees; this is the content of the next result.

\begin{cor}
\label{cor:gowers3}
Let $\F \subseteq W^{<\omega}(\Sigma ; \upsilon)$ be a family of finite sequences of variable words 
over an alphabet $\Sigma$. Then for every infinite sequence $\vec{w} \in W^\omega (\Sigma ; \upsilon)$ of variable words over $\Sigma$ there exists a variable reduction $\vec{u}\prec\vec{w}$ of $\vec{w}$ over $\Sigma$ such that:

 either $VRW^{<\omega}(\vec{u})\subseteq \F$, 

 or for every variable reduction $\vec{u}_1 $ of $\vec{u}$ there exists an initial segment of $\vec{u}_1$ 
which belongs to $W^{<\omega}(\Sigma ; \upsilon)\setminus \F$. 
\end{cor}
\begin{proof} Let $\F_1 = \{\bt=(t_1,\ldots,t_k) \in \F: (t_1,\ldots,t_n)\in\F$ for every $n\in\nat$ with $n\leq k \}\cup \{\emptyset\}$.
\newline
The family $\F_1$ is a tree. According to Corollary~\ref{cor:gowers}, there exists $\vec{u}\prec\vec{w}$ such that:

 either $VR^{<\omega}(\vec{u})\subseteq \F_1 \subseteq \F$, 

 or for every variable reduction $\vec{u}_1 $ of $\vec{u}$ there exists an initial segment of $\vec{u}_1$ 
which belongs to $\W^{<\omega}(\Sigma ; \upsilon)\setminus \F_1$. 
\newline
Let $\vec{u}_1=(u^1_n)_{n\in\nat}\prec \vec{u}$, and let $k\in\nat$ such that $\bt=(u^1_1,\ldots,u^1_k)\in W^{<\omega}(\Sigma ; \upsilon)\setminus \F_1 = (\F\setminus \F_1)\cup (W^{<\omega}(\Sigma ; \upsilon)\setminus \F)$. 
Then, either $\bt \in W^{<\omega}(\Sigma ; \upsilon)\setminus \F$, as required, 
or $\bt\in \F\setminus \F_1$. 
In case $\bt\in \F\setminus \F_1$, by the definition of $\F_1$, there 
exists $n\in\nat$ with $n\leq k$ such that $(u^1_1,\ldots,u^1_n)\in W^{<\omega}(\Sigma ; \upsilon) \setminus\F$,
as required.
\end{proof}

The result for families of (constant) words, corresponding to Corollary~\ref{cor:gowers3}, can now be obtained as a corollary to Theorem~\ref{block-NashWilliams1}. 
\begin{cor}
\label{cor:gowers1} 
Let $\G\subseteq W^{<\omega}(\Sigma)$ be a family of finite sequences of words over the alphabet $\Sigma$. Then for every infinite sequence $\vec{w} \in W^\omega (\Sigma ; \upsilon)$ of variable words over $\Sigma$ there exists a variable reduction $\vec{u}\prec\vec{w}$ of $\vec{w}$ over $\Sigma$ such that:

 either $RW^{<\omega}(\vec{u})\subseteq \G$, 

 or for every reduction $\vec{u}_1 $ of $\vec{u}$ there exists an initial segment of $\vec{u}_1$ 
which belongs to $W^{<\omega}(\Sigma)\setminus \G$. 
\end{cor}
\begin{proof}
If $\G$ is a tree, then the proof follows from Theorem~\ref{block-NashWilliams1} (case 1 and subcase 2(ii)) and Proposition~\ref{prop:canonicalrep2}, if not set 

$\G_1 = \{\bs=(s_1,\ldots,s_l) \in \G: (s_1,\ldots,s_k)\in\G$ for every $k\in\nat$ with $k\leq l \}\cup \{\emptyset\}$.
\newline
The family $\G_1$ is a tree and $\G_1 \subseteq \G$. Hence, there exists $\vec{u}\prec\vec{w}$ such that:

 either $RW^{<\omega}(\vec{u})\subseteq \G_1\subseteq \G$, 

 or for every reduction $\vec{u}_1 $ of $\vec{u}$ there exists an initial segment of $\vec{u}_1$ 
which belongs to $W^{<\omega}(\Sigma)\setminus \G_1$. 
\newline
Given that for every reduction $\vec{u}_1 $ of $\vec{u}$ there exists an initial segment $\bs$ of $\vec{u}_1$ 
which belongs to $W^{<\omega}(\Sigma)\setminus \G_1=(W^{<\omega}(\Sigma)\setminus \G)\cup(\G\setminus \G_1)$ we have that there exists an initial segment of of $\vec{u}_1$ which belongs to $W^{<\omega}(\Sigma) \setminus\G$,
as required.
\end{proof}

Corollary~\ref{cor:gowers1} is equivalent to the following infinitary partition 
theorem, which is the counterpart for (constant) words of Carlson's infinitary partition theorem (Theorem~\ref{cor:IblockNW}, Corollary~\ref{cor:gowers}).

\begin{cor} 
\label{cor:blockNW1}
Let $\U \subseteq W^{\omega}(\Sigma)$ be a pointwise closed family of infinite sequences of words over a finite non-empty alphabet $\Sigma$ and $\vec{w}\in W^\omega (\Sigma ; \upsilon)$ be an infinite sequence of variable words.  Then there exists a variable reduction $\vec{u}\prec\vec{w}$ of $\vec{w}$ over $\Sigma$ such that:
$$\text{either }\ RW^{\omega}(\vec{u})\subseteq \U\ \quad \text{or}\quad 
RW^{\omega}(\vec{u})\subseteq W^{\omega}(\Sigma) \setminus \U\ .$$
\end{cor}

\section{Schreier-type version of Carlson's Ellentuck type partition theorem for words}
In this final section, we establish (in Theorem~\ref{thm:Ellentuck}) a rather technical strengthening of Theorem B (mentioned in the introduction) derived from Theorem~\ref{block-NashWilliams2}, involving the Ellentuck topology $\Tau_E$, defined on $W^{\omega}(\Sigma ; \upsilon)$ (Definition~\ref{def:Ellentuck}). A simple consequence of Theorem~\ref{thm:Ellentuck} is the characterization of completely Ramsey partitions of $W^{\omega}(\Sigma ; \upsilon)$ in terms of the Baire property in the topology $\Tau_E$ (Corollary~\ref{cor:Baireproperty}), a result proved with different methods by Carlson in \cite{C}. A similar characterization of completely Ramsey partitions of $W^{\omega}(\Sigma)$ can be proved analogously, as a consequence of Theorem~\ref{block-NashWilliams1}.  

We start by defining the topology $\Tau_E$ on $W^{\omega}(\Sigma ; \upsilon)$, an 
analogue of the Ellentuck topology on $[\nat]^\omega$, defined in \cite{E}. For simplicity, we write 
$\emptyset \propto \vec{w}$ and $\vec{w} \setminus \emptyset = \vec{w}$ for every $\vec{w}\in W^{\omega}(\Sigma ; \upsilon)$. 

\begin{defn}\label{def:Ellentuck}
Let $\Tau_E$ be the topology on $W^{\omega}(\Sigma ; \upsilon)$ with basic open sets of 
the form $[\bs ,\vec{s}]$ for $\bs \in W^{<\omega}(\Sigma ; \upsilon)$ and 
$\vec{s}\in W^{\omega}(\Sigma ; \upsilon)$, where for $\bs\in W^{<\omega}(\Sigma ; \upsilon) \setminus \{\emptyset\}$
\begin{center}
$[\bs, \vec{s}] =  \{\vec{w} \in W^{\omega}(\Sigma ; \upsilon): \bs \propto \vec{w}\quad\text{and}
\quad \vec{w} \setminus \bs \prec \vec{s}\}$ and $[\emptyset, \vec{s}] =  \{\vec{w} \in W^{\omega}(\Sigma ; \upsilon):  \vec{w} \prec \vec{s}\}$.
\end{center}
The topology $\Tau_E$ is stronger than the relative topology of 
$W^{\omega}(\Sigma ; \upsilon)$ with respect to the pointwise convergence topology of 
$\{0,1\}^{[D]^{<\omega}}$, which has basic open sets of the form 
$[\bs, \vec{e}] = \{\vec{w} \in W^{\omega}(\Sigma ; \upsilon): \bs\propto \vec{w}\}$ for $\bs \in W^{<\omega}(\Sigma ; \upsilon)$ and $\vec{e}=(e_n)_{n \in \nat}$ with $e_n=\upsilon$ for every $n\in\nat$.

We denote by $\hat\U$ and $\U^\lozenge$ the closure and the interior 
respectively of a family $\U\subseteq W^{\omega}(\Sigma ; \upsilon)$ in the topology 
$\Tau_E$. Then it is easy to see that 
\begin{equation*}
\begin{split}
&\hat\U = \{\vec{w}\in W^{\omega}(\Sigma ; \upsilon): [\bs, \vec{w} \setminus \bs]\cap \U\ne\emptyset
\ \text{ for every }\  \bs\propto \vec{w}\}\ ;\ \text{ and}\\
&\U^\lozenge = \{\vec{w}\in W^{\omega}(\Sigma ; \upsilon):\ \text{ there exists }
\bs\propto \vec{w}\ \text{ such that }\ 
[\bs, \vec{w} \setminus \bs] \subseteq \U\}\ .
\end{split}
\end{equation*}
\end{defn}

Now we can state the main theorem of this section. For $ \bs = (s_1,\ldots,s_k), 
\bt = (t_1,\ldots,t_l) \in W^{<\omega}(\Sigma ; \upsilon)$  we set $ \bs \odot \bt = (s_1,\ldots,s_k,t_1,\ldots,t_k)$  and $ \emptyset \odot \bt = \bt \odot \emptyset=\bt $.  

\begin{thm}\label{thm:Ellentuck}
Let $\U\subseteq W^{\omega}(\Sigma ; \upsilon)$, $\bs\in W^{<\omega}(\Sigma ; \upsilon)$ and 
$\vec{s}\in W^{\omega}(\Sigma ; \upsilon)$. 
Then

 either there exists $\vec{u} \prec \vec{s}$ such that $[\bs,\vec{u}]\subseteq \hat\U$,

 or there exists a countable ordinal $\xi_0= \zeta_{(\bs,\vec{s})}^{\U}$ 
such that for every $\xi>\xi_0$ there exists $\vec{u} \prec \vec{s}$ with 
$[\bs \odot \bt,\vec{t} \setminus \bt]\subseteq W^{\omega}(\Sigma ; \upsilon)\setminus \U$ for every 
$\vec{t} \prec \vec{u}$ and $\bt \in W^\xi(\Sigma ; \upsilon)$ with $\bt \propto \vec{t}$.
\end{thm}

We will give the proof of this theorem, using Lemma~\ref{lem:Ellentuck1}, an analogue of Lemma~\ref{lem:Carlson-Ramsey 2}. 

\begin{lem}\label{lem:Ellentuck1}
Let $\R \subseteq\{[\bs,\vec{s}] :\bs\in W^{<\omega}(\Sigma ; \upsilon), \vec{s}\in W^{\omega}(\Sigma ; \upsilon)\}$ with the properties:

{\rm (i)} for every $(\bs,\vec{s})\in W^{<\omega}(\Sigma ; \upsilon)\times W^{\omega}(\Sigma ; \upsilon)$
there exists $\vec{s}_1 \prec \vec{s}$ such that $[\bs,\vec{s}_1]\in\R$; and 

{\rm (ii)} for every $[\bs,\vec{s}_1]\in\R$ and $\vec{s}_2 \prec \vec{s}_1$ we have 
$[\bs,\vec{s}_2] \in\R$.
\newline
Then, for every $(\bs,\vec{s})\in \W^{<\omega}(\Sigma ; \upsilon)\times\W^{\omega}(\Sigma ; \upsilon)$ 
there exists $\vec{u}\in [\bs, \vec{s}]$ such that $[\bs\odot\bt,\vec{t} \setminus \bt]\in\R$ 
for every $\vec{t} \prec \vec{u} \setminus \bs$ and $\bt \propto \vec{t}$.
\end{lem}

\begin{proof}
Let $ \bs \in \W^{<\omega}(\Sigma ; \upsilon)$ and $\vec{s}\in\W^{\omega}(\Sigma ; \upsilon)$. 
According to assumption~(i), there exists $\vec{s}_1=(s_n^1)_{n\in\nat} \prec \vec{s}$ such that 
$[\bs, \vec{s}_1] \in\R$. Set $u_1=s_1^1\in VRW(\vec{s})$. According to the assumption~(i), there exists 
$\vec{s}_2=(s_n^2)_{n\in\nat} \prec \vec{s}_1\setminus (u_1)$ such that $[\bs \odot (u_1), \vec{s}_2]\in \R$. Set 
$u_2=s_1^2$. Then $(u_1,u_2)\in VRW^{<\omega}(\vec{s})$ and $\vec{s}_2  \prec \vec{s}\setminus u_1$. 

Let $n\in \nat$, $n>1$. Assume that there have been constructed $\vec{s}_1, \ldots, \vec{s}_n \in\W^{\omega}(\Sigma ; \upsilon)$ and $u_1,u_2,\ldots,u_n \in W(\Sigma ; \upsilon)$ such that $(u_1,u_2,\ldots, u_n)\in VRW^{<\omega}(\vec{s})$, $\vec{s}_{i+1} \prec \vec{s}\setminus u_1\ast\ldots\ast u_{i}$, 
$\vec{s}_{i+1} \prec \vec{s}_{i} \setminus (u_{i})$ for every $i=1,\ldots,n-1$, $(u_i) \propto \vec{s}_i$ for every $i=1,\ldots,n$, and 
$[\bs \odot \bt, \vec{s}_{i+1}]\in\R$ for every $\bt \in VRW^{<\omega}((u_1,\ldots, u_i))$ and $i=1,\ldots,n-1$.
 
We will construct $\vec{s}_{n+1}$ and $u_{n+1}$. Let $\{\bt_1,\ldots,\bt_m\} = VRW^{<\omega}((u_1,\ldots, u_n))$ for some $m\in \nat$. According to assumption~(i), there exist 
$\vec{s}_{n+1}^1,\ldots, \vec{s}_{n+1}^m \in W^\omega (\Sigma ; \upsilon)$ such that 
$\vec{s}_{n+1}^m \prec \cdots \prec \vec{s}_{n+1}^1 \prec \vec{s}_n \setminus (u_n)$ and 
$[\bs \odot \bt_i, \vec{s}_{n+1}^i]\in \R$ for every $i=1,\ldots,m$. Set $\vec{s}_{n+1} = \vec{s}_{n+1}^m$. If $\vec{s}_{n+1} = (s^{n+1}_i)_{i\in\nat}$, then set $u_{n+1} = s^{n+1}_1$. Of course $\vec{s}_{n+1}\prec \vec{s}_n \setminus (u_n)$, $(u_{n+1}) \propto \vec{s}_{n+1}$, $(u_0,u_1,\ldots, u_{n+1})\in VRW^{<\omega}(\vec{w})$, 
$\vec{s}_{n+1} \prec \vec{s}\setminus u_1\ast\ldots\ast u_{n}$ and, according to condition (ii), 
$[\bs \odot \bt, \vec{s}_{n+1}]\in\R$ for every $\bt \in VRW^{<\omega}((u_1,\ldots, u_n))$.

Set $\vec{u} = (s_1,\ldots,s_k,u_1,u_2,\ldots)\in\W^{\omega}(\Sigma ; \upsilon)$ in case $\bs=(s_1,\ldots,s_k)$ and 
$\vec{u} = (u_1,u_2,\ldots)$ in case $\bs=\emptyset$. 
Then $\vec{u} \in [\bs,\vec{s}]$, since $(u_1,u_2,\ldots, u_n)\in VRW^{<\omega}(\vec{s})$ for every $n\in\nat$. 
Let $\vec{t} \prec \vec{u} \setminus \bs$ and $\bt\ne\emptyset$, $\bt \propto \vec{t}$. Since 
$\bt\in VRW^{<\omega}(\vec{u} \setminus \bs)$, set $n_0 = \min \{n\in\nat: \bt \in VRW^{<\omega}((u_1,\ldots,u_n))\}$.
Then $[\bs \odot \bt,\vec{s}_{n_0+1}]\in \R$. According to assumption~(ii), $[\bs \odot \bt, \vec{u} \setminus (s_1,\ldots,s_k,u_1,\ldots,u_{n_0})]\in \R$ and 
$[\bs \odot \bt,\vec{t}\setminus\bt]\in\R$, since 
$\vec{t}\setminus\bt \prec \vec{u} \setminus \bs \odot (u_1,\ldots,u_{n_0}) \prec \vec{s}_{n_0+1}$. If $\bt=\emptyset$, then $[\bs,\vec{u}\setminus \bs] \in\R$, since $[\bs,\vec{s}_1] \in\R$ and, according to assumption~(ii),  $[\bs,\vec{t}] \in\R$.
\end{proof}

\begin{proof}[Proof of Theorem~\ref{thm:Ellentuck}] 
Set 
\begin{equation*}
\begin{split}
\R_{\U} = &\{[\bt,\vec{t}] :(\bt,\vec{t}) \in W^{<\omega}(\Sigma ; \upsilon)\times W^{\omega}(\Sigma ; \upsilon)
\ \text{ and}\\
&\qquad \text{either }\ [\bt,\vec{t}]\cap \U = \emptyset\quad\text{or}\quad 
[\bt,\vec{t}_1]\cap \U\ne\emptyset\ \text{ for every }\ \vec{t}_1 \prec \vec{t}\}\ .
\end{split}
\end{equation*}
It is easy to check that $\R_{\U}$ satisfies the assumptions (i) and (ii) 
of Lemma~\ref{lem:Ellentuck1}, 
hence, there exists $\vec{u}=(u_n)_{n\in\nat}\in [\bs, \vec{s}]$ such that $[\bs\odot\bt,\vec{t} \setminus \bt]\in\R_{\U}$ for every $\vec{t} \prec \vec{u} \setminus \bs$ and $\bt \propto \vec{t}$.

For $\bt\in VRW^{<\omega}(\vec{u} \setminus \bs)$ there exists unique $\vec{u}_\bt \prec \vec{u} \setminus \bs$ with $\bt \propto \vec{u}_\bt$ and $\vec{u}_\bt \setminus \bt=(u_n)_{n>n_0}$ for some $n_0\in\nat$. Then 
$[\bs\odot\bt,\vec{u}_\bt \setminus \bt]\in\R_{\U}$. Set

$\F = \{\bt\in VRW^{<\omega}(\vec{u} \setminus \bs) : [\bs \odot \bt, \vec{u}_1] \cap \U \ne \emptyset\text{ for every } \vec{u}_1 \prec \vec{u}_\bt \setminus \bt\}$.

The family $\F$ is a tree. 
Indeed, let $\bt \in \F$ and $\bt_1 \propto \bt$.
Then $[\bs \odot \bt_1, \vec{u}_{\bt_1} \setminus \bt_1] \in\R_{\U}$, since 
$\bt_1 \in VRW^{<\omega}(\vec{u} \setminus \bs)$. 
It is impossible $[\bs \odot \bt_1,\vec{u}_{\bt_1} \setminus \bt_1]\cap \U = \emptyset$, since 
$[\bs \odot \bt,\vec{u}_1]\cap \U \ne \emptyset$ for every $\vec{u}_1 \prec \vec{u}_\bt \setminus \bt$. Hence, 
$[\bs \odot \bt_1, \vec{u}_{1}]\cap \U \ne \emptyset$ for every $\vec{u}_1 \prec \vec{u}_{\bt_1} \setminus \bt_1$, and consequently $\bt_1 \in \F$. We now apply Theorem~\ref{block-NashWilliams2} for $\F$ and 
$\vec{u} \setminus \bs$. We have the following cases:

[Case 1] There exists $\vec{u}_1 \prec \vec{u} \setminus \bs \prec \vec{s}$ such that 
$VRW^{<\omega}(\vec{u}_1) \subseteq \F$. 
This gives that $[\bs \odot \bt, \vec{u}_2]\cap \U \ne \emptyset$ for every 
$\bt \in VRW^{<\omega}(\vec{u}_1)$ and $\vec{u}_2 \prec \vec{u}_1\setminus \bt$, 
which implies that $[\bs, \vec{u}_1] \subseteq \hat\U$.

[Case 2] There exists a countable ordinal $\xi_0 = \zeta_{\vec{u} \setminus \bs}^{\F} 
= \zeta_{(\bs,\vec{s})}^{\U}$ such that for every $\xi>\xi_0$ there exists 
$\vec{u}_1 \prec \vec{u} \setminus \bs \prec \vec{s}$ with 
$W^{\xi}(\Sigma ; \upsilon)\cap VRW^{<\omega}(\vec{u}_1) \subseteq W^{<\omega}(\Sigma ; \upsilon)\setminus \F$. Using the canonical representation of every infinite sequense of variable words with respect to the family $W^{\xi}(\Sigma ; \upsilon)$ (Proposition~\ref{prop:canonicalrep}), we have that $[\bs \odot \bt, \vec{u}_1] \subseteq W^{\omega}(\Sigma ; \upsilon)\setminus \U$ for every 
$\bt\in W^{\xi}(\Sigma ; \upsilon)\cap VRW^{<\omega}(\vec{u}_1)$,. Hence, $[\bs \odot \bt, \vec{t} \setminus \bt]\subseteq W^{\omega}(\Sigma ; \upsilon)\setminus \U$ for every 
$\vec{t} \prec \vec{u}_1$ and $\bt \in W^\xi(\Sigma ; \upsilon)$ with $\bt \propto \vec{t}$.
\end{proof}

Applying Theorem~\ref{thm:Ellentuck} to partitions $\U$ that are closed 
or meager in the topology $\Tau_E$, we obtain the following consequences. 

\begin{cor}\label{cor:Ellentuck2}
Let $\U$ be a closed in the topology $\Tau_E$ subset of $W^{\omega}(\Sigma ; \upsilon)$, $\bs\in W^{<\omega}(\Sigma ; \upsilon)$ and $\vec{s}\in W^{\omega}(\Sigma ; \upsilon)$. Then

 either there exists $\vec{u} \prec \vec{s}$ such that $[\bs,\vec{u}]\subseteq \U$,

 or there exists a countable ordinal $\xi_0$ 
such that for every $\xi>\xi_0$ there exists $\vec{u} \prec \vec{s}$ with 
$[\bs \odot \bt,\vec{t} \setminus \bt]\subseteq W^{\omega}(\Sigma ; \upsilon)\setminus \U$ for every 
$\vec{t} \prec \vec{u}$ and $\bt \in W^\xi(\Sigma ; \upsilon)$ with $\bt \propto \vec{t}$.
\end{cor}

\begin{cor}\label{cor:Ellentuck3}
Let $\U$ be a subset  of $W^{\omega}(\Sigma ; \upsilon)$ meager in the topology $\Tau_E$, 
$\bs \in W^{<\omega}(\Sigma ; \upsilon)$ and $\vec{s}\in W^{\omega}(\Sigma ; \upsilon)$. 
Then, there exists a countable ordinal $\xi_0$ such that for every $\xi>\xi_0$
there exists $\vec{u} \prec \vec{s}$ with 
$[\bs \odot \bt,\vec{t} \setminus \bt]\subseteq W^{\omega}(\Sigma ; \upsilon)\setminus \U$ for every 
$\vec{t} \prec \vec{u}$ and $\bt \in W^\xi(\Sigma ; \upsilon)$ with $\bt \propto \vec{t}$.
\end{cor}

\begin{proof}
We use Theorem~\ref{thm:Ellentuck} for $\U$.
We will prove that the first alternative is impossible. 
Indeed, let $\vec{u} \prec \vec{s}$ such that $[\bs, \vec{u}]\subseteq \hat\U$.  
If $\U= \bigcup_{n\in\nat} \U_n$ with $(\hat\U_n)^\lozenge = \emptyset$ 
for every $n\in\nat$, then we set 
\begin{equation*}
\begin{split}
\R = &\{[\bt,\vec{t}] : \bt \in W^{<\omega}(\Sigma ; \upsilon),\ \vec{t}\in W^{\omega}(\Sigma ; \upsilon)
\ \text{ and}\\
&\qquad [\bt,\vec{t}]\cap \U_k =\emptyset\ \text{ for every } \ k\in\nat
\text{ with } k\le |\bt|\};\ 
\end{split}
\end{equation*}
where $|\bt|$ denotes the number of terms of the finite sequence $\bt\in W^{<\omega}(\Sigma ; \upsilon)$. The family $\R$ obviously satisfies the assumption (ii) of Lemma~\ref{lem:Ellentuck1} and also satisfies the assumption (i) of Lemma~\ref{lem:Ellentuck1}, since, according to Theorem~\ref{thm:Ellentuck} and Proposition~\ref{prop:canonicalrep}, for every $\bt\in W^{<\omega}(\Sigma ; \upsilon)$, $\vec{t}\in W^{\omega}(\Sigma ; \upsilon)$ and $k\in\nat$ there exists $\vec{t}_1 \prec \vec{t}$ such that $[\bt, \vec{t}_1]\cap \U_k = \emptyset$, as 
it is impossible $[\bt, \vec{t}_1]\subseteq \hat\U_k$ for every $k\in\nat$. Hence, according to Lemma~\ref{lem:Ellentuck1}, there exists 
$\vec{u}_1 \in [\bs,\vec{u}]$ such that 
$[\bs \odot \bt,\vec{u}_1 \setminus \bs \odot \bt]\in\R$ for every $\bt\in VRW^{<\omega}(\vec{u}_1 \setminus \bs)$. 

We will prove that $[\bs, \vec{u}_1 \setminus \bs]\cap \U=\emptyset$. 
Let $\vec{u}_2\in [\bs,\vec{u}_1 \setminus \bs]\cap\U$. 
Then, $\vec{u}_2 \in [\bs, \vec{u}_1 \setminus \bs] \cap \U_k$ for some $k\in\nat$. 
Hence, there exists $\bt\in VRW^{<\omega}(\vec{u}_1 \setminus \bs)$ such that $\bs \odot \bt\propto \vec{u}_2$, 
$k\le |\bs \odot \bt|$ and $[\bs \odot \bt, \vec{u}_1 \setminus \bs \odot \bt]\cap \U_k\ne\emptyset$.
Then, $[\bs \odot \bt,\vec{u}_1 \setminus \bs \odot \bt]\notin \R$. 
A contradiction, since $\bt\in VRW^{<\omega}(\vec{u}_1 \setminus \bs)$. 
Hence, $[\bs,\vec{u}_1 \setminus \bs]\cap \U=\emptyset$, and consequently $\vec{u}_1\notin \hat\U$.

This is a contradiction, since $\vec{u}_1 \in [\bs ,\vec{u}] \subseteq\hat\U$. 
Hence, the first alternative of Theorem~\ref{thm:Ellentuck} for the partition $\U$ is 
impossible, so the second alternative holds for $\U$.
\end{proof}

We recall the definition of the completely Ramsey families of infinite sequences of variable words given by Carlson in \cite{C}.

\begin{defn}\label{def:Milliken}
A family $\U\subseteq W^{\omega}(\Sigma ; \upsilon)$ of infinite sequences of variable words on a finite, non-empty  alphabet $\Sigma$ is called {\em completely Ramsey} if for every $\bs\in W^{<\omega}(\Sigma ; \upsilon)$ and every $\vec{s}\in W^{\omega}(\Sigma ; \upsilon)$ there exists $\vec{u} \prec \vec{s}$ such that 
\begin{equation*}
\text{either }\ [\bs,\vec{u}]\subseteq \U\quad\text{ or }\quad 
[\bs,\vec{u}]\subseteq W^{\omega}(\Sigma ; \upsilon) \setminus \U\ .
\end{equation*}
\end{defn}

The characterization of completely Ramsey families of infinite sequences of variable words, proved with different methods by Carlson in \cite{C}, is a consequence of Theorem~\ref{thm:Ellentuck}.

\begin{cor} [\textsf{Carlson}, \cite{C}]
\label{cor:Baireproperty}
A family $\U\subseteq W^{\omega}(\Sigma ; \upsilon)$ is completely Ramsey if and only if 
$\U$ has the Baire property in the topology $\Tau_E$.
\end{cor} 

\begin{proof} 
Let $\U\subseteq W^{\omega}(\Sigma ; \upsilon)$ have the Baire property in the topology 
$\Tau_E$. 
Then $\U = \B \triangle  \C  = (\B \cup \C^c) \cup (\C\cap \B^c)$, where 
$\B\subseteq W^{\omega}(\Sigma ; \upsilon)$ is $\Tau_E$-closed and 
$\C\subseteq W^{\omega}(\Sigma ; \upsilon)$ is $\Tau_E$-meager $(\C^c= W^{\omega}(\Sigma ; \upsilon) \setminus \C)$.  
According to Corollary~\ref{cor:Ellentuck3} and Proposition~\ref{prop:canonicalrep}, 
for every $\bs\in W^{<\omega}(\Sigma ; \upsilon)$ and $\vec{s}\in W^{\omega}(\Sigma ; \upsilon)$, 
there exists $\vec{u} \prec \vec{s}$ such that $[\bs,\vec{u}]\subseteq \C^c$ and according to Corollary~\ref{cor:Ellentuck2}, there exists $\vec{u}_1 \prec \vec{u}$ such that 
\begin{itemize}
\item[{}] either $[\bs,\vec{u}_1] \subseteq \B\cap [\bs,\vec{u}]\subseteq \B\cap \C^c \subseteq \U$;
\item[{}] or $[\bs,\vec{u}_1]\subseteq \B^c \cap [\bs,\vec{s}] \subseteq \B^c \cap \C^c \subseteq \U^c$.
\end{itemize}
Hence, $\U$ is completely Ramsey.

On the other hand, if $\U$ is completely Ramsey, then $\U = \U^\lozenge
\cup (\U\setminus \U^\lozenge)$ and $\U\setminus\U^\lozenge$ is a meager set 
in $\Tau_E$. Hence $\U$ has the Baire property in the topology $\Tau_E$.
\end{proof}

A similar characterization of the completely Ramsey families of $W^{\omega}(\Sigma)$ can be proved analogously, as a consequence of Theorem~\ref{block-NashWilliams1}.

\begin{remark}\label{rem:compx}
(i) The Ellentuck topology $\Tau_E$ on $W^{\omega}(\Sigma)$ has basic open sets of the form $[\bs ,\vec{s}]$ for $\bs \in W^{<\omega}(\Sigma )$ and $\vec{s}\in W^{\omega}(\Sigma ; \upsilon)$, where   
\begin{center}
$[\bs, \vec{s}] =  \{\vec{w} \in W^{\omega}(\Sigma): \bs \propto \vec{w}\quad\text{and}
\quad \vec{w}\setminus \bs \in RW^{\omega}(\vec{s})\}$ and $[\emptyset, \vec{s}] = RW^{\omega}(\vec{s})$. 
\end{center}

(ii) A family $\U\subseteq W^{\omega}(\Sigma)$ of infinite sequences of words on a finite, non-empty alphabet $\Sigma$ is called {\em completely Ramsey} if for every $\bs\in W^{<\omega}(\Sigma)$ and every $\vec{s}\in W^{\omega}(\Sigma ; \upsilon)$ there exists $\vec{u} \prec \vec{s}$ such that 
\begin{equation*}
\text{either }\ [\bs,\vec{u}]\subseteq \U\quad\text{ or }\quad 
[\bs,\vec{u}]\subseteq W^{\omega}(\Sigma) \setminus \U\ .
\end{equation*}

(iii) A family $\U\subseteq W^{\omega}(\Sigma)$ is completely Ramsey if and only if 
$\U$ has the Baire property in the topology $\Tau_E$.
\end{remark}

\bigskip
{\footnotesize
\noindent 
{\sc Department of Mathematics, Athens University, Panepistemiopolis, Athens 15784, Greece}
\noindent
\newline  
E-mail address of first named author: vfarmaki@math.uoa.gr
\newline 
E-mail address of second named author: snegrep@math.uoa.gr
}

\end{document}